\documentclass[english,12pt]{article}

\ifx\optionkeymacros\undefined\else \fi

\catcode`\Œ=\active\defŒ{{\aa}}       
\catcode`\º=\active\defº{\int}        
\catcode`\=\active\def{\c c}        
\catcode`\¶=\active\def¶{\partial}    
\catcode`\Ä=\active\defÄ{\oint}       
\catcode`\Æ=\active\defÆ{\triangle}   
\catcode`\Â=\active\defÂ{\neg}        
\catcode`\µ=\active\defµ{\mu}         
\catcode`\¿=\active\def¿{{\o}}        
\catcode`\¹=\active\def¹{\pi}         
\catcode`\Ï=\active\defÏ{{\oe}}       
\catcode`\§=\active\def§{{\ss}}       
\catcode`\ =\active\def {\dagger}     
\catcode`\Ã=\active\defÃ{\sqrt}       
\catcode`\·=\active\def·{\Sigma}      
\catcode`\Å=\active\defÅ{\approx}     
\catcode`\½=\active\def½{\Omega}      
\catcode`\£=\active\def£{{\it\$}}     
\catcode`\°=\active\def°{\infty}      
\catcode`\¤=\active\def¤{{\S}}        
\catcode`\¦=\active\def¦{{\P}}        
\catcode`\¥=\active\def¥{\bullet}     
\catcode`\»=\active\def»{\leavevmode\raise.585ex\hbox{\b a}}      
\catcode`\¼=\active\def¼{\leavevmode\raise.6ex\hbox{\b o}}        
\catcode`\­=\active\def­{\not=}       
\catcode`\²=\active\def²{\leq}        
\catcode`\³=\active\def³{\geq}        
\catcode`\Ö=\active\defÖ{\div}        
\catcode`\É=\active\defÉ{{\dots}}     
\catcode`\¾=\active\def¾{{\ae}}       
\catcode`\Ç=\active\defÇ{\ll}         
\catcode`\Ò=\active\defÒ{``}          
\catcode`\Á=\active\defÁ{!`}          
\catcode`\¢=\active\def¢{\rlap/c}     
\catcode`\Ô=\active\defÔ{`}           
\catcode`\Õ=\active\defÕ{'}           

\catcode`\=\active\def{{\AA}}       
\catcode`\'=\active\def'{\c C}        
\catcode`\¯=\active\def¯{{\O}}        
\catcode`\¸=\active\def¸{\Pi}         
\catcode`\Î=\active\defÎ{{\OE}}       
\catcode`\®=\active\def®{{\AE}}       
\catcode`\×=\active\def×{\diamond}    
\catcode`\¡=\active\def¡{\accent'27}  
\catcode`\Ó=\active\defÓ{''}          
\catcode`\±=\active\def±{\pm}         
\catcode`\È=\active\defÈ{\gg}         
\catcode`\À=\active\defÀ{?`}          
\catcode`\Ð=\active\defÐ{--}          
\catcode`\Ñ=\active\defÑ{---}         

\catcode`\Š=\active\defŠ{\"a}        
\catcode`\'=\active\def'{\"e}        
\catcode`\•=\active\def•{\"{\i}}     
\catcode`\š=\active\defš{\"o}        
\catcode`\Ÿ=\active\defŸ{\"u}        
\catcode`\Ø=\active\defØ{\"y}        
\catcode`\€=\active\def€{\"A}        
\catcode`\…=\active\def…{\"O}        
\catcode`\†=\active\def†{\"U}        
\catcode`\‡=\active\def‡{\'a}        
\catcode`\Ž=\active\defŽ{\'e}        
\catcode`\'=\active\def'{\'{\i}}     
\catcode`\—=\active\def—{\'o}        
\catcode`\œ=\active\defœ{\'u}        
\catcode`\ƒ=\active\defƒ{\'E}        
\catcode`\ˆ=\active\defˆ{\`a}        
\catcode`\=\active\def{\`e}        
\catcode`\"=\active\def"{\`{\i}}     
\catcode`\˜=\active\def˜{\`o}        
\catcode`\=\active\def{\`u}        
\catcode`\Ë=\active\defË{\`A}        
\catcode`\‹=\active\def‹{\~a}        
\catcode`\–=\active\def–{\~n}        
\catcode`\›=\active\def›{\~o}        
\catcode`\Ì=\active\defÌ{\~A}        
\catcode`\"=\active\def"{\~N}        
\catcode`\Í=\active\defÍ{\~O}        
\catcode`\‰=\active\def‰{\^a}        
\catcode`\=\active\def{\^e}        
\catcode`\"=\active\def"{\^{\i}}     
\catcode`\™=\active\def™{\^o}        
\catcode`\ž=\active\defž{\^u}        

\let\optionkeymacros\null

\usepackage[english]{babel}
\input{amssym.def}
\input{amssym.tex}

\newcommand{\binomial}[2]{\displaystyle{{#1 \choose #2}}}

\newcommand{\sqm}[4]{\displaystyle{
\left({#1 \atop #3}{#2 \atop #4}\right)}}

\def\AA{{\Bbb A}}
\def\RR{{\Bbb R}}
\def\CC{{\Bbb C}}
\def\QQ{{\Bbb Q}}
\def\NN{{\Bbb N}}
\def\PP{{\Bbb P}}
\def\ZZ{{\Bbb Z}}
\def\GG{{\Bbb G}}
\def\FF{{\Bbb F}}
\def\TT{{\Bbb T}}
\def\tr{{\bf t}}
\def\no{{\bf n}}

\newcommand{\GL}{\mbox{\rm GL}}

\newcommand{\Ker}{\mbox{\rm Ker}}

\newcommand{\und}{\underline}
\newcommand{\spt}{\Sigma}

\newcounter{constant}

\newcounter{fonction}

\newtheorem{Definition}{Definition}[section]
\newtheorem{Lemme}{Lemma}[section]
\newtheorem{Proposition}{Proposition}

\begin{document}
\begin{center}
\selectlanguage{english}

\LARGE{On the arithmetic properties of complex\\ values of Hecke-Mahler series.}\\

\vspace{15pt}

\Large{Federico Pellarin}
\end{center}

\vspace{10pt}

\noindent{\footnotesize {\bf Abstract.} 
Here we characterise in a complete and explicit way
the relations of algebraic dependence over $\QQ$ of complex values of
Hecke-Mahler series taken at algebraic points of the multiplicative group
$\GG_m^2(\CC)$. Our result contains previous theorems by 
Loxton and  van der Poorten, Mahler, and Masser.}

\begin{footnotesize}
\tableofcontents
\end{footnotesize}


\newpage

\section{Introduction, results.\label{section:introduction}}

Let $w$ be a real positive irrational number. The {\em Hecke-Mahler
series associated to} $w$ is the power series:
\[f_w(u,v)=\sum_{l=1}^\infty\sum_{h=1}^{[lw]}u^lv^h,\]
where the square brackets denote the greatest integer part.
This series, which is transcendental, converges for all $u,v$ 
complex numbers such that $|u|<1$ and $|u||v|^w<1$.

On p. 208 of \cite{Masser:Hecke}, Masser asked several questions which may
all be included in the following

\vspace{10pt}

\noindent {\bf Problem.} {\em Given an $m$-tuple ${\cal M}=((u_1,v_1),
\ldots,(u_m,v_m))=(\und{u}_1,\ldots,\und{u}_m)$ 
of non-zero complex numbers, compute the transcendence degree over the field of rational numbers
$\QQ$, of the following subfield of the field of complex numbers:
\[\QQ(\und{u}_1,\ldots,\und{u}_m,f_{w_1}(\und{u}_1),\ldots,f_{w_m}(\und{u}_m)),\]
for positive irrationals $w_i$ (when all the complex numbers above make sense).}

\vspace{10pt}

The aim of this text is to completely solve the particular aspect of the problem above, which 
consists of choosing $w_1=\cdots=w_m=w$ a {\em quadratic} irrational, and 
the couples of complex numbers $(u_i,v_i)$ {\em algebraic} over $\QQ$.

More precisely,
we introduce a certain explicit condition of geometric nature,
on an $m$-tuple of couples of non-zero algebraic numbers 
\begin{equation}
{\cal M}=((u_1,v_1),\ldots,(u_m,v_m)),\label{eq:mtuple}\end{equation}
that we call Òsemi-freeness": this condition depends on $w$ and will be made explicit
in special cases only, to simplify our exposition. 

In theorem 2 we prove, under some technical hypotheses on $w$, that an $m$-tuple
${\cal M}$ as above is semi-free if and only if the complex 
numbers $f_w(u_1,v_1),\ldots,f_w(u_m,v_m)$ are algebraically independent over $\QQ$.

These technical hypotheses on $w$ are harmless: in theorem 3 of the appendix, we 
completely characterise all the $m$-tuples $((u_1,v_1),\ldots,(u_m,v_m))$ as above,
such that the complex 
numbers $f_w(u_1,v_1),\ldots,f_w(u_m,v_m)$ are algebraically independent over $\QQ$
for a general quadratic irrational $w>0$.

Our results improve a theorem of Masser in \cite{Masser:Hecke}:

\vspace{10pt}

\noindent {\bf Theorem 1 (Masser).}
{\em Let $u_1,\ldots,u_m$ be algebraic numbers such that $0<|u_i|<1$ for all $i=1,\ldots,m$.
Then $f_w(u_1,1),\ldots,f_w(u_m,1)$ are algebraically independent if and only if
$u_1,\ldots,u_m$ are distinct.}

\vspace{10pt}

In \cite{Masser:Hecke}, Masser predicted the relations between complex numbers such as $f_w(u,v)$ to be
Òquite complicated": here are some examples of relations.

If we take $m=5$ and $(u,v)$ a couple of
complex numbers such that $|u|<1$ and $0<|u||v|^w<1$, then the following homogeneous linear
relation can be easily checked:
\begin{equation}
4f_w(u^2,v^2)-f_w(u,v)-f_w(-u,v)-f_w(u,-v)-f_w(-u,-v)=0,\label{eq:masser}
\end{equation}
This relation is in some sense as simple as possible, because it holds for 
any choice of $u,v,w$. 

For certain quadratic irrationals $w$ only, there also exists 
Òshorter" relations. Indeed, there sometimes also exist
positive rational integers $a,b,c,d$, with
$\det\sqm{a}{b}{c}{d}=1$, and a rational function $R(u,v)\in\QQ(u,v)$, such that:
\begin{equation}
f_w(u,v)-f_w(u^av^b,u^cv^d)=R(u,v).
\label{eq:fu_eq}\end{equation} For this special choice of $w$, this is
a {\em functional equation} of $f_w$:
thus, if $u,v$ are algebraic, then we get non-homogeneous linear relations, this time with $m=2$.

In section \ref{section:portrait} of the appendix,
other Òspecial" and Ògeneric" relations will be explicitely described, 
and we can make right now a commentary about them. 

We will find that all algebraic relations are generated by linear relations.
moreover, these linear relations are always connected with homogeneous
linear forms of rational functions.
The so-called {\em Hecke's geometric series} (as defined in \cite{Hecke:Analytische}) are
also helpful to encode the relations (see section \ref{section:linear}).

But the easiest way to describe the relations is to employ 
formal series of a certain type, with rational
coefficients, that will be introduced later; this suggests our definition of semi-freeness (definition
\ref{defi:semi_free}). 

If we compare the relations
that we find with, for example, the relations that conjecturally connect the complex values of the logarithm (Shanuel's
conjecture), generated by:
\[\log(uv)-\log(u)-\log(v)=0,\quad u,v\in\CC^\times\] then definitely, the relations connecting the special complex
values of Hecke-Mahler series are quite complicated. By the way, 
there is some analogy between some of these relations and the relations connecting special values
of Fredholm series, as in \cite{Loxton:Fredholm}.

Here are a few words about the methods employed in this article. 
The main difficulty we have to overcome is the fact that we
must work with  analytic functions of two complex variables (theorem 1 essentially deals with
functions of one complex variable).

Mahler's method, when it applies, is an eccellent technique to investigate arithmetic
properties of functions of several variables: that is the way we attack our problem. Thus, the functional
equation (\ref{eq:fu_eq}) will have a priviledged meaning, and this explains why we 
need $w$ to be quadratic, satisfying certain hypotheses.

In our proof of theorem 2, we use a classical criterion of algebraic independence of Loxton and 
van der Poorten; this criterion is used in many articles about the arithmetic properties of
values of locally analytic functions (\footnote{The notion of locally analytic function will be given in section
\ref{section:semifree}.}) satisfying certain functional equations. For example,
the same criterion is used in \cite{Masser:Hecke}.
We also use the powerful vanishing theorem of Masser as in \cite{Masser:Vanishing}
to check a technical condition,
called Òproperty A", playing the role of a zero estimate. The vanishing theorem is also used in \cite{Masser:Hecke}.
We do not need to generalise these tools in order to obtain our results; in particular
we do not need to perform any particular construction of transcendence here;
arithmetically, the classical tools of Mahler's theory are enough to prove 
our results.

The new point in our approach is to use certain structures of real multiplications on tori.
In particular we interpret the functional equation (\ref{eq:fu_eq}) by using an algebraic action of a unit of $K$ on a
product of multiplicative groups. All the relations between complex values of $f_w$ will be nicely 
described by using algebraic actions of orders of $K$ on tori.

The introduction in the
theory of these new tools leads to new problems. Analytically and algebraically,
these actions deserve
quite a few surprises, as the reader will see.
We hope that ours will be a good viewpoint also for similar or more general problems.

We now introduce the condition of semi-freeness; the main theorem will be stated in section \ref{section:main_theorem}.

\subsection{The condition of Òsemi-freeness".}

Let $\bar{\QQ}$ denote the algebraic closure of $\QQ$, let us fix an
embedding of $\bar{\QQ}$ in the field of complex numbers $\CC$.
Let $K\subset\bar{\QQ}$ be a real quadratic number field.
The chosen embedding $\bar{\QQ}\rightarrow\CC$ induces an embedding 
$\sigma:K\rightarrow\RR$: in all the following we will consider
$K$ as a subfield of $\RR$; for example, the expressions $\nu>0$ and $\nu'>0$ mean, 
for $\nu\in K$, $\sigma(\nu)>0$ and $\sigma(\nu')>0$ respectively.

If $\nu$ is an element of $K$ we denote
$\nu'$ its non-trivial Galois conjugate. From now on, the field $K$ is fixed.
It should be said right now, that one of the most important elementary facts used in this article
is that if $\eta>1$ is a unit of $K$, then:
\[\lim_{k\rightarrow\infty}\eta'{}^k=0.\]

We denote by $\GG_m(\CC)=\CC^\times$ the complex multiplicative group.
In this article, we will work in the group:
\[\TT:=\GG_m^2(\CC),\] 
with identity element $\und{1}=(1,1)$,
and in its powers $\TT^n=\TT\oplus\cdots\oplus\TT$, for some $n\in\NN$.

All throughout this text the elements of $\CC^n,\RR^n,\ldots$ are considered as
row matrices, unless otherwise specified. If ${\cal A}_1,\ldots,{\cal A}_n$ are square matrices we 
denote ${\cal A}_1\oplus\cdots\oplus{\cal A}_n$ the square matrix having diagonal blocks equal
to the ${\cal A}_i$'s, and zero elsewhere.

\subsubsection{Exponential functions.}

Let $M$ be a complete $\ZZ$-module of $K$, that is, a free 
$\ZZ$-module of rank $2$ contained in $K$, let us fix a $\ZZ$-basis $(B_0,B_{1})$
of $M$. We note $M^*$ the dual of $M$ for the trace $\tr:K\rightarrow\QQ$ that is, 
the complete $\ZZ$-module  of $K$ whose elements $\nu$ satisfy $\tr(\mu\nu)\in\ZZ$ 
for all $\mu\in M$; if $M\subset N$ are complete $\ZZ$-modules, then
$N^*\subset M^*$. 

We will sometimes extend the map $\tr:K\rightarrow\QQ\subset\CC$, as follows.
If we have a	couple of complex numbers $\und{z}=(z,z')$, we will often write:
\[\tr(\und{z})=z+z'.\] The underlined expression $\und{z}$ will be often simplified, and we will most of the time
write $z$ instead of $\und{z}$. This, of course, makes no sense, but nicely 
simplifies plenty of expressions. If $\nu\in K$, for example, the expression
$\tr(\nu z)$ will mean $\nu z+\nu' z'$.

We note $(B_0^*,B_1^*)$ the dual basis of $(B_0,B_{1})$
for the trace and we note:
\begin{equation}{\goth B}=\sqm{B_0^*}{B_0^*{}'}{B_1^*}{B_1^*{}'}=\sqm{B_0}{B_1}{B_0'}{B_1'}^{-1}.
\label{eq:goth_B}\end{equation} Let $\Sigma:K\rightarrow\RR^2$
be the embedding $\Sigma(\nu)=(\sigma(\nu),\sigma(\nu'))\in\RR^2\subset\CC^2$.
We have the exponential function with periods in $\Sigma(M)$
\[\Phi:\CC^2\rightarrow\TT,\]
defined by~:
\begin{equation}
\Phi(z,z')={}^{{\rm t}}e({\goth B}\cdot{}^{{\rm t}}(z,z')),
\label{eq:exponentielle}\end{equation}
where ${}^{{\rm t}}\cdot$ means Òtranspose", and $e(\tau)=e(2\pi{\rm i}\tau)$ (with ${\rm i}:=\sqrt{-1}$). 
In a more explicit way, if
$(u,v)=\Phi(z,z')$, then: 
\begin{equation}
u=e( B_0^*z+B^*{}'_0z'),\quad v=e( B_1^*z+B^*{}_1'z').\label{eq:more_explicit}\end{equation}
This function $\Phi$ factors through $\CC^2/\Sigma(M)$
because for complex numbers $z,z',\zeta,\zeta'$ we have ${\goth B}\cdot{}^{{\rm t}}(z,z')-{\goth B}\cdot{}^{{\rm
t}}(\zeta,\zeta')\in\ZZ^2$ if and only if $(z-\zeta,z'-\zeta')\in\Sigma(M)$, and this happens if and only if
$\Phi(z,z')=\Phi(\zeta,\zeta')$. Of course, if $(z,z'),(\tilde{z},\tilde{z}')\in\CC^2$, then we have, in $\TT$:
\[\Phi(z,z')\Phi(\tilde{z},\tilde{z}')=\Phi(z+\tilde{z},z'+\tilde{z}').\]

\subsubsection{Actions of orders of $K$ on $\TT$.}

Let \[{\cal B}=\sqm{a}{b}{c}{d}\in\GL_2(\QQ)\] be a regular matrix 
with rational integer entries $a,b,c,d$,
let $\und{u}=(u,v)$ be an element of $\TT(\CC)$.
We denote: \[{\cal B}.\und{u}=(u^av^b,u^cv^d).\] 

Let $S=S(M)$ be the stabiliser of
$M$, i. e. the order of
$K$ whose elements are the $\beta$'s such that $\beta M\subset M$. Notice that if $\nu\in K\setminus\{0\}$, then 
$S(\nu M)=S(M)$.

The multiplicative group $\TT$ is endowed with an action of 
$S$ analytically defined as follows: let $\und{u}=(u,v)$ be a point of $\TT$,
let $(z,z')\in\CC^2$ be such that $\Phi(z,z')=(u,v)$ and
let $\mu$ be an element of $S$. Then:
\begin{equation}
\und{u}^\mu:=\Phi(\mu z,\mu'z').\label{eq:action_S}\end{equation} This action is well defined and
depends on $M$ as well as on the basis $(B_0,B_1)$; moreover, the action is algebraic. Indeed, 
let $\mu$ be an element of $S-\{0\}$; let us denote:
\begin{equation}
{\cal B}(\mu)={\goth B}\cdot\sqm{\mu}{0}{0}{\mu'}\cdot{\goth B}^{-1}.\label{eq:explicit_1}
\end{equation}
A simple computation shows that ${\cal B}(\mu)$ has determinant $\no(\mu)$ 
(where $\no(\mu)$ is the norm $\mu\mu'$ of $\mu$ over $\QQ$) and 
has rational integer coefficients. Moreover, we have that \[\und{u}^\mu={\cal B}(\mu).\und{u}.\]

With $M$ fixed, we have several actions (\ref{eq:action_S}) which might behave in very different ways
(we will give some examples later in
the section
\ref{section:example} of the appendix).
Nevertheless, if $\mu\in\ZZ$, then 
\begin{equation}\und{u}^\mu=(u^\mu,v^\mu)\label{eq:Z_action}\end{equation} does not even depend on $K$,
and equals the usual componentwise $\mu$-power in $\TT$ (the usual diagonal action of $\ZZ$ on $\TT$).

\begin{Definition}{\em 
We say that a point $\und{u}\in\TT$ is a {\em torsion point} if $\und{u}=(\zeta_1,\zeta_2)$
with $\zeta_1,\zeta_2$ roots of unit. A point which is not a torsion point is said to be a
{\em point of infinite order.}}\end{Definition}

It is easy to see that $\und{u}$ is torsion if and only if: \begin{equation}\und{u}\in\Phi(\Sigma(K)),
\label{eq:torsion_sigma}\end{equation} that is, there exists
$\alpha\in K$ such that $\und{u}=\Phi(\alpha,\alpha')$. The subgroup of $\TT$ whose elements are torsion points
is denoted by $\TT_{{\tiny \mbox{tors}}}$.

\subsubsection{Formal series.}
Let
$\und{U}=(U,V)$ be a couple of indeterminates. Formal power series 
\[\sum_{{\tiny \begin{array}{c}
h,l\in\ZZ\\ (h,l)\not=(0,0)\end{array}}}c_{h,l}U^hV^l\] with $c_{h,l}\in\CC$ live in a $\CC$-vector 
space in the usual way, regardless to their convergence. 

We must introduce a certain class of formal power series. We first observe that the power series:
\[\FF(\und{U})=\sum_{{\tiny \begin{array}{c}
h,l\in\ZZ\\ (h,l)\not=(0,0)\end{array}}}U^hV^l\] only converges at $U=V=0$; in the following we will consider
this (and other similar series) as a formal power series in the space defined above. Notice that
we do not identify it with a rational function.

Let $\nu$ be an element of $M^*$. There exist $h,l\in\ZZ$ (unique) such that
\[\nu=hB_0^*+lB_1^*.\]
Let us write $M(\und{U})^\nu=U^hV^l$; we then have \[\FF(\und{U})=\sum_{\nu\in M^*\setminus\{0\}}M(\und{U})^\nu.\]
We now consider, more generally, a complete $\ZZ$-module $N$ of $K$, containing $M$ (so that $N^*\subset M^*$) and we attach
to it the formal power series:
\[\FF_N(\und{U})=\sum_{\nu\in N^*\setminus\{0\}}M(\und{U})^\nu.\] As a special case, we obtain $\FF_M=\FF$.
\subsubsection{Definition of semi-freeness.}
We fix an exponential function $\Phi$
as in (\ref{eq:exponentielle}), and we remark that, if $\und{u},\und{v}\in\TT$ have infinite order, then only two cases can
occur.
\begin{enumerate}
\item 
We have
\begin{equation}\und{u}^\gamma=\Phi(\alpha,\alpha')\und{v}^\beta,\label{eq:buono}\end{equation}
for some
$\alpha\in K$ and
$\beta,\gamma\in S\setminus\{0\}$.
\item For all $\beta,\gamma\in S$, not all zero, $\und{u}^\beta\und{v}^{\gamma}$ has infinite order.
\end{enumerate}
Let us consider an $m$-tuple ${\cal M}$ as in (\ref{eq:mtuple}); the relation $\approx$ determines a 
partition $({\cal J}_k)_{k=1,\ldots,s}$ of the set $\{1,\ldots,m\}$.

There exists $\und{v}^\sharp\in\TT$, such that $\und{v}^\sharp{}^\gamma=\und{v}$, with 
$\gamma$ as in (\ref{eq:buono}). Thus, we can find $\delta\in K,\rho\in S\setminus\{0\}$ such that 
$\und{u}=\Phi(\delta,\delta')\und{v}^\sharp{}^\rho$ and (\ref{eq:buono}) holds.

We define an equivalence relation $\approx$ on the subset of elements of infinite
order of $\TT$ by declaring that $\und{u}\approx\und{v}$ if we are in case 1. 

A little induction, and the argument above, show that for all $k=1,\ldots,s$, there exists $\und{v}_k\in\TT$ of infinite
order,  such that for all $j\in{\cal J}_k$, there exists $\alpha_j\in K$ and $\beta_j\in S\setminus\{0\}$ 
with the property that:
\begin{equation}
\und{u}_j=\Phi(\alpha_j,\alpha'_j)\und{v}_k^{\beta_j}\label{eq:torsion_and_not}
\end{equation}

We associate to $\und{u}_j$ a 
formal power series $\FF_{\und{v}_k}(\Phi,\und{u}_j:\und{U})$ depending on the choice of $\und{v}_k$ and $\Phi$.
If the relation (\ref{eq:torsion_and_not}) holds, then we set:
\[\FF_{\und{v}_k}(\Phi,\und{u}_j:\und{U})=
\FF_{\beta_j^{-1}M}\left(\Phi\left(\frac{\alpha_j}{\beta_j},\frac{\alpha_j'}{\beta_j'}\right)\und{U}\right).\]
\begin{Definition}{\em Let \[{\cal M}=(\und{u}_1,\ldots,\und{u}_m)\in\TT^m\] be a $m$-tuple
of elements of infinite order of $\TT$.
We say that ${\cal M}$ is {\em semi-free} with respect to $\Phi$, if 
for all $k=1,\ldots,s$, the series:
\[\FF_{\und{v}_k}(\Phi,\und{u}_j:\und{U}),\quad j\in{\cal J}_k\]
are $\QQ$-linearly independent.}
\label{defi:semi_free}\end{Definition}

If an $m$-tuple ${\cal M}=(\und{u}_1,\ldots,\und{u}_m)$ as above
is semi-free with respect to $\Phi$, we will sometimes say that the 
elements $\und{u}_1,\ldots,\und{u}_m$ are {\em multiplicatively semi-independent} with respect to $\Phi$.

\medskip

\noindent {\bf Remark.} The series $\FF_{\und{v}_k}(\Phi,\und{u}_j:\und{U})$ depends on the 
choice of $\und{v}_k$; but the property of semi-freeness does not. This will be proved in lemma \ref{lemme:two_groups} of the
appendix.

\subsection{A special choice $\Phi_0$ for $\Phi$ and the  main theorem.\label{section:main_theorem}}

To state our theorem we have to associate to the quadratic
irrational $w$ a complete $\ZZ$-module $M$ and an exponential function $\Phi_0$ with
periods in it. 

\vspace{10pt}

We state our theorem for $w$'s with particular properties. In section \ref{section:general_w}
of the appendix, we will show
that with these restrictive hypotheses, we can also deal the case of general 
quadratic irrationals $w$.

Let us suppose that $w=\theta$ is such that:
\begin{equation}
0<\theta<1\mbox{ and }\theta'<-1.
\label{eq:perron}\end{equation}
Let $M=\ZZ+\theta^{-1}\ZZ$ be the complete $\ZZ$-module with the basis $(B_0,B_1)=(\theta^{-1},1)$,
let ${\goth B}_0$ be the matrix ${\goth B}$ as in (\ref{eq:goth_B}), for the basis above:
\[{\goth B}_0=\sqm{B_0^*}{B_0^*{}'}{B_1^*}{B_1^*{}'}=
\Delta^{-1/2}\sqm{1}{-1}{-\theta'{}^{-1}}{\theta^{-1}},\] where $\Delta$ is the
discriminant of $M$, let us consider the exponential function:
\begin{equation}
\Phi_0(z,z')={}^{{\rm t}}e( {\goth B}_0\cdot{}^{{\rm t}}(z,z')).
\label{eq:exponential_function}\end{equation}

\vspace{10pt}

\noindent {\bf Theorem 2.} {\em Let
${\cal M}=((u_1,v_1),\ldots,(u_m,v_m))$ be an $m$-tuple
of algebraic elements of $\TT$ such that
$|u_i|<1$ and $0<|u_i||v_i|^\theta<1$ for $i=1,\ldots,m$. Then ${\cal M}$ is 
semi-free with respect to $\Phi_0$ if and only if the complex numbers 
$f_\theta(u_1,v_1),\ldots,f_\theta(u_m,v_m)$ 
are algebraically independent over $\QQ$. 

Moreover, if $f_\theta(u_1,v_1),\ldots,f_\theta(u_m,v_m)$ 
are algebraically dependent over $\QQ$ then there is
a non trivial linear relation \begin{equation}\sum_{i=1}^mc_if_\theta(u_i,v_i)=\lambda,\label{eq:relations_1}\end{equation} where
$\lambda$ is an algebraic number and $c_1,\ldots,c_m$ are rational numbers.}

\vspace{10pt}

In section \ref{section:portrait2}, we will describe in detail the relations
(\ref{eq:relations_1}); explicit examples will be given in section \ref{section:portrait} of the appendix. 
Some examples about the different behavior
of two distinct exponential functions  
will be given in section \ref{section:example} of the appendix.

The condition of 
semi-freeness (definition \ref{defi:semi_free}) is in general rather difficult to
check, but this is quite natural because non-trivial linear relations with many terms can occur between the numbers
$f_\theta(u_i,v_i)$; one example is provided by (\ref{eq:masser}). However, there
are some natural additional conditions on ${\cal M}$ which have the property
that they discard all the non-trivial relations with strictly more than two terms $f_\theta(u_i,v_i)$.

\vspace{10pt}

\noindent {\bf Corollary 1.} {\em Let $H$ be a connected algebraic subgroup of $\TT$ of 
dimension $1$. Let $(u_1,v_1),\ldots,(u_m,v_m)$ be algebraic elements of $H$
such that $|u_i|<1$ and $0<|u_i||v_i|^\theta<1$ for $i=1,\ldots,m$.
Then $f_\theta(u_1,v_1),\ldots,f_\theta(u_m,v_m)$ are algebraically independent over $\QQ$ if
and only if $(u_1,v_1),\ldots,(u_m,v_m)$ are distinct.}

\vspace{10pt}

The corollary 1 implies theorem 1 by considering $H=\GG_m(\CC)\times\{1\}$, and will be proved in section
\ref{section:consequences} of the appendix.

\vspace{10pt}

\noindent {\bf Corollary 2.} {\em Let $(u_1,v_1),\ldots,(u_m,v_m)$ be algebraic elements of $\TT$
such that $|u_i|<1$ and $0<|u_i||v_i|^\theta<1$ for $i=1,\ldots,m$. Let us suppose that, via the exponential function
$\Phi_0$, the elements $(u_1,v_1),\ldots,(u_m,v_m)$ 
generate an $S$-submodule of $\TT$ which is
contained  in an $S$-submodule of $\TT$ isomorphic to 
a finite direct sum $S\oplus\cdots\oplus S$. Then we have that
$(u_i,v_i)^\eta\not=(u_j,v_j)$ for all $1\leq i<j\leq m$ and for all units $\eta$ of $S$, if and only if 
$f_\theta(u_1,v_1),\ldots,f_\theta(u_m,v_m)$ are
algebraically independent over $\QQ$.}

\vspace{10pt}

 The corollary 2 implies a result of Loxton and van der Poorten
in \cite{Loxton:Arithmetic} (see theorem 3) and a proof will also be given in section
\ref{section:consequences}. 

\medskip

We fix right away a quadratic irrational $\theta\in K$ satisfying 
(\ref{eq:perron}), the complete $\ZZ$-module $M=\ZZ+\theta^{-1}\ZZ$ with the basis 
\[(B_0,B_1)=(\theta^{-1},1),\] and the exponential function $\Phi_0$ as in (\ref{eq:exponential_function}). 
Since $\theta$ is fixed all along this text (except in section \ref{section:general_w} of the appendix), we denote
\[f(\und{u}):=f_\theta(\und{u}).\] 
We denote ${\cal D}$ the domain of convergence of $f(\und{u})$; that is:
\[{\cal D}=\{(u,v)\in\CC^2\mbox{ such that }|u|<1\mbox{ and }|u||v|^\theta<1\}.\]
Together with $f$, we need to consider the twin series:
\begin{eqnarray*}
f^+(\und{u})& =& \sum_{l=1}^\infty\sum_{h=[\theta l]+1}^{[(2/\tr(\theta^{-1}))l]}u^lv^h.
\end{eqnarray*}
It is easy to check that the series $f^+$ converges in the domain:
\[{\cal D}^+=\{(u,v)\in\CC^2\mbox{ such that $|u|<1$ and $|u|^{\tr(\theta^{-1})}|v|^{2}<1$}\}.\]

\section{Algebraic tools.\label{section:semifree}}
The proof of theorem 2 remains on the algebraic study of quite a few functions; here are the associated symbols:
\[f_N,\quad f_N^+,\quad\Theta_N,\quad R_{\alpha,\beta,N},\quad R_{\alpha,\beta,N}^+.\] 
We are going to define them, and to study their rationality properties.

\medskip

Let \begin{equation}F(Z,Z')=\sum_{\nu\in M^*}c_\nu e(\tr(\nu Z)),
\label{eq:F1}\end{equation} be a formal series, let
$\und{U}=(U,V)$ be a couple of variables
in $\TT$ formally connected to the couple of variables $(Z,Z')$ by $\und{U}=\Phi_0(Z,Z')$. To $F$ is associated 
a formal power series $G$ un the variables $U,V$, defined by:
\[G(\und{U})=F(Z,Z').\] The formal series $G$ is clearly well defined, because the periods of $\Phi_0$ lie
in $\Sigma(M)$, and if $(z,z'),(\zeta,\zeta')$ differ by an element of $\Sigma(M)$, then for all
$\nu\in M^*$, $e(\tr(\nu z))=e(\tr(\nu \zeta))$.
\begin{Definition}{\em Let $F$ as in (\ref{eq:F1}). The {\em $K$-support} $\spt_K(F)$
is the subset of $\Sigma(M^*)\subset\RR^2$ whose
elements are the $\Sigma(\nu)$'s such that $c_\nu\not=0$. Let $G$ be associated to $F$ as above.
The $K$-support $\spt_K(G)$ of $G$ is equal, by definition, to the $K$-support of $F$.}
\label{def:def_K_spt}\end{Definition}

Here are more notations.
\begin{eqnarray*}
K_+&=&\{\nu\in K\mbox{ such that }\nu>0\mbox{ and }\nu'>0\},\\
K_\pm &=&\{\nu\in K\mbox{ such that }\nu>0\mbox{ and
}\nu'<0\},\\
{\cal I}&=&\{\nu\in K\mbox{ such that }\nu>-\nu'>0\mbox{ or
}\nu\geq\nu'>0\}.
\end{eqnarray*}
If $E$ is any subset of $K$, we denote:
\begin{eqnarray*}
E_+&=&E\cap K_+,\\
E_\pm&=&E\cap K_\pm.
\end{eqnarray*} 
Thus, ${\cal I}={\cal I}_\pm\cup{\cal I}_+$, with ${\cal I}_\pm=\{\nu\in K\mbox{ such that }\nu>-\nu'>0\}$
and ${\cal I}_+=\{\nu\in K\mbox{ such that }\nu\geq\nu'>0\}$.

We now introduce new series which somewhat generalise the series $f$. 
Let $N$ be a complete $\ZZ$-module of $K$ containing $M$, let $\und{u}=\Phi_0(z,z')$, let us define
the series:
\begin{eqnarray}
f_N(\und{u}) & = & \sum_{\nu\in{\cal I}_\pm\cap N^*}e( \tr(\nu z)),\label{eq:no1}\\
f^+_N(\und{u}) & = & \sum_{\nu\in{\cal I}_+\cap N^*}e( \tr(\nu z)),\label{eq:no2}\\
\Theta_N(\und{u}) & = & \sum_{\nu\in{\cal I}\cap N^*}e( \tr(\nu z))=f_N(\und{u})+f_N^+(\und{u}).\label{eq:theta_1}
\end{eqnarray}
These series are clearly well defined, thanks to the discussion above;
they are uniquely determined by their $K$-supports. In figure 1, these $K$-supports are represented
together and compared.

Let ${\cal H}$ be the complex upper half plane
\[{\cal H}=\{z\in\CC\mbox{ such that }\Im(z)>0\},\]
let us introduce two subsets of $\CC^2$:
\begin{eqnarray*}
{\cal W} & = & \{(z,z')\in{\cal H}\times\CC\mbox{ with }\Im(z')<\Im(z)\},\\
{\cal W}^+ & = & \{(z,z')\in{\cal H}\times\CC\mbox{ with }-\Im(z)<\Im(z')<\Im(z)\},
\end{eqnarray*} 
where $\Im(\cdot)$ denotes the imaginary part 
of a complex number; clearly,
\[{\cal W}^+\subset{\cal W}.\]

It is easy to check that, for all $N$, the series $f_N$ converges for $\und{u}\in\TT$ such that
$\und{u}=\Phi_0(z,z')$, with $(z,z')\in{\cal W}$, and
the series $f_N^+,\Theta_N$ both converges for $\und{u}$ such that $(z,z')\in{\cal W}^+$.

\medskip

So far, we have introduced three types of series $f_N,f_N^+$ and $\Theta_N$; we need two more types.

Let $\alpha\in K, \beta\in S-\{0\}$ such that $\beta>0$, 
let $N$ be a complete $\ZZ$-module containing $M$, such that $S(N)=S$, and let us denote:
\[\und{\zeta}=\Phi_0(\alpha,\alpha'),\quad\und{\zeta}^\sharp=\Phi_0(\alpha/\beta,\alpha'/\beta').\] If $\beta\in S_+$ we
set:
\begin{eqnarray*}
R_{\alpha,\beta,N}(\und{u})&=&f_{\beta^{-1}N}(\und{\zeta}^\sharp\und{u})-f_N(\und{\zeta}\und{u}^\beta)\\
R_{\alpha,\beta,N}^+(\und{u})&=&f_{\beta^{-1}N}^+(\und{\zeta}^\sharp\und{u})-f^+_N(\und{\zeta}\und{u}^\beta).
\end{eqnarray*}
Otherwise $\beta\in S_\pm$, and we set:
\begin{eqnarray*}
R_{\alpha,\beta,N}(\und{u})&=&
f_{\beta^{-1}N}(\und{\zeta}^\sharp\und{u})-f^+_N(\und{\zeta}\und{u}^\beta)\\
R_{\alpha,\beta,N}^+(\und{u})&=&
f_{\beta^{-1}N}^+(\und{\zeta}^\sharp\und{u})-f_N(\und{\zeta}\und{u}^\beta).\end{eqnarray*}
To simplify our notations, we will write $R_{\alpha,\beta},R_{\alpha,\beta}^+$ instead of
$R_{\alpha,\beta,M},R_{\alpha,\beta,M}^+$. Moreover, if $\alpha=0$, we will drop the corresponding subscript
so that we will more simply write $R_\beta$ instead of $R_{0,\beta}$ and $R_\beta^+$ instead of $R_{0,\beta}^+$;
similarly, $R_{\beta,N}=R_{0,\beta,N}$ and $R_{\beta,N}^+=R_{0,\beta,N}^+$.

The main series of this article having been defined, we may now study their basic property in the next section.

\subsection{Rational locally analytic functions.}

\begin{Definition}{\em A {\em locally analytic function} on $\CC^{n}$ is a function
defined over a subset of $\CC^{n}$ which 
is analytic on a non-empty open neighbourhood of the origin
$\und{0}=(0,\ldots,0)\in \CC^{n}$.}\end{Definition}

We start with an elementary lemma.
\begin{Lemme}
Let $\rho_1\in\RR$ and $\rho_2\in\RR_{\geq 0}\cup\{\infty\}$ be such that 
\begin{equation}-\rho_2<\rho_1<\rho_2,\label{eq:conditions_1}\end{equation} 
let $Z\subset\ZZ^2$ be a subgroup of finite index.
Then, the series
\begin{equation}R(u,v)=\sum_{{\tiny \begin{array}{c}(l,h)\end{array}}}u^lv^h,
\label{eq:series_rho}\end{equation} where the sum runs over the
set 
\[\{(l,h)\in Z\mbox{ such that }l>0\mbox{ and }\rho_1l<h\leq\rho_2l\}\] 
is locally analytic.

If $\rho_1,\rho_2\in\PP_1(\QQ)$ (that is, if $\rho_1\in\QQ$ and $\rho_2\in\QQ_{\geq 0}\cup\{\infty\}$), then
$R(u,v)\in\QQ(u,v)$.
\label{lemme:rational_functions}\end{Lemme}
\noindent {\bf Proof.} 
The convergence of the series $R$ is easily checked:
its domain of convergence is the set \[\{(u,v)\in\CC^2\mbox{ such that $|u||v|^{\min\{0,\rho_1\}}<1$ and
$|u||v|^{\rho_2}<1$}\}.\] This is clearly a non-empty open neighbourhood of $\und{0}\in\CC^2$ 
because of condition (\ref{eq:conditions_1}) on $\rho_1,\rho_2$. We continue with the rationality property.

Let us first suppose that $Z=\ZZ^2$. 
If $(a,b),(c,d)$ are couples of coprime rational integers such that
$\rho_1=a/b, \rho_2=c/d$ and $a,c\geq 0$, then:
\[(1-u^av^b)(1-u^cv^d)\sum_{l=1}^\infty\sum_{h=[\rho_1l]+1}^{[\rho_2l]}u^lv^h\]
is a finite sum $\sum_{i,j}u^iv^j$, where $(i,j)\in\ZZ^2$, which proves that $R$ is rational.

In the case where $Z$ has finite index in $\ZZ^2$,
we observe that there exists a matrix $A\in\GL_2(\QQ)$ with
rational integer coefficients such that $R(u,v)=R'(A.(u,v))$ where $R'$ is a 
sum like (\ref{eq:series_rho}), with $Z=\ZZ^2$, thus 
lying in $\QQ(u,v)$, which proves that $R$ is a rational function in this case
too.

\begin{Proposition}
Let $N$ be a complete $\ZZ$-module containing $M$, and satisfying $S(N)=S$; we have the following properties.
\begin{enumerate}
\item The series $f_N,f_N^+,\Theta_N$ define locally analytic functions in $\CC^2$. 
\item If $N=M$, then $f_M=f$ and
$f_M^+=f^+$.
\item The series $\Theta_N$ defines a rational function in $\QQ(\und{u})$.
\item The $K$-support of $R_{\alpha,\beta,N}$ is contained in $M^*_\pm$, and the $K$-support 
of $R_{\alpha,\beta,N}^+$ is contained in $M^*_+$.
\item For all $\alpha\in K$, $\beta\in S$ with $\beta>|\beta'|>0$, the series 
$R_{\alpha,\beta,N},R_{\alpha,\beta,N}^+$ define
locally analytic rational functions in $\bar{\QQ}(\und{u})$.
\end{enumerate}
\label{lemme:theta_etc}\end{Proposition}
\noindent {\bf Proof.} 
As a preliminary remark, we note that if $(u,v)=\Phi_0(z,z')$, then:
\begin{equation}
u^lv^h=e( \tr(\nu z)),
\label{eq:usone}\end{equation} where $\nu\in M^*$ is defined by
\begin{equation}\nu=\Delta^{-1/2}(-h\theta'{}^{-1}+l)=hB_1^*+lB_0^*.\label{eq:nu_delta}\end{equation}

\medskip

\noindent {\bf (1)}. Lemma \ref{lemme:rational_functions} allows to see that the series 
$f_N,f_N^+,\Theta_N$ define locally analytic functions in $\CC^2$. Their domains of convergence can be determined 
explicitly; the domain of convergence of the series $f_N$ is the set ${\cal D}$, whereas the 
common domain of convergence of the series $f_N^+,\Theta_N$ is the set ${\cal D}^+$. 

\medskip

\noindent {\bf (2)}.
We prove that $f_M=f$. 
We have that $\nu\in{\cal I}_\pm\cap M^*$ if and only if $\nu>-\nu'>0$ and $h,l\in\ZZ$, which is equivalent,
using (\ref{eq:usone}), to:
\[\Delta^{-1/2}h(\theta^{-1}-\theta'{}^{-1})>0\mbox{ and }l>\theta^{-1}h,\]
that is, taking into account that $\rho(1)=0$,
\begin{equation}
h>0,\quad h<\theta l.
\label{eq:tonto}\end{equation}
Since $\theta^{-1}-\theta'{}^{-1}>0$ (condition (\ref{eq:perron})), $\nu\in{\cal I}_\pm\cap M^*$ if and only if 
$(h,l)\in\ZZ^2$ and $l\geq 1\mbox{ and }0< h<\theta l\}$, and this ensures that $f_M=f$.

The proof of the equality $f_M^+=f^+$ is very similar,
and we omit it. 

\medskip

\noindent {\bf (3)}. Let us first suppose that $N=M$. Then, by definition, $\Theta_M=f+f^+$, so that:
\[\Theta_M(\und{u})=\sum_{l=1}^\infty\sum_{h=1}^{[(2/\tr(\theta^{-1}))l]}u^lv^h.\]
Applying lemma \ref{lemme:rational_functions} with $\rho_1=0,\rho_2=2/\tr(\theta^{-1})\in\QQ_{>0}$,
we see that $\Theta_M\in\QQ(\und{u})$. For general $N$ the proof is the same because to $N$ is associated 
a $\ZZ$-submodule of finite index $Z\subset \ZZ^2$, and lemma \ref{lemme:rational_functions} can be applied to it.

\medskip

\noindent {\bf (4)}. Let $\beta\in S$, such that $\beta>0$. Then, $\beta\in S_+$ or $\beta\in S_\pm$.
If $\beta\in S_+$, then $f_N(\und{\zeta}\und{u}^\beta)$ has its $K$-support entirely contained in
$M^*_\pm$, and this, for every torsion point $\und{\zeta}$. The same happens to $f_N^+(\und{\zeta}\und{u}^\beta)$ 
when $\beta\in
S_\pm$. Thus, for all $\beta>0$ and for all $\alpha$, $R_{\alpha,\beta,N}$ has its support in $M^*_\pm$, because it
is a linear combination of series already having this property.
The proof of the property for the series $R_{\alpha,\beta,N}^+$ is similar.

\medskip

\noindent {\bf (5)}. 
Let $\nu\in K-\{0\}$, let us denote:
\[\rho(\nu)=\frac{\tr(\nu\Delta^{-1/2})}{\tr(\theta^{-1}\nu\Delta^{-1/2})}\in\QQ\cup\{\infty\}
\mbox{ and }\rho^+(\nu)=\frac{\tr(\nu)}{\tr(\theta^{-1}\nu)}\in\QQ\cup\{\infty\}.\]
By using (\ref{eq:perron}), it is straightforward to see that
\begin{eqnarray}
\rho^+(\beta)>\theta>\rho(\beta)>0, & & \mbox{ if }\beta>\beta'>0,\label{eq:111}\\
\rho(\beta)>\theta>\rho^+(\beta)>0, & & \mbox{ if }\beta>-\beta'>0\label{eq:222}.
\end{eqnarray}
Let us suppose that $\beta>\beta'>0$. With the notations introduced above:
\begin{eqnarray}
\lefteqn{R_{\alpha,\beta,N}(\Phi_0(z,z'))=}\nonumber\\ & = & \sum_{\nu\in \beta N^*\cap{\cal
I}_\pm}e\left(\tr\left(\nu\left(\frac{\alpha}{\beta}+z\right)\right)\right)-
\sum_{\nu\in N^*\cap{\cal I}_\pm}e(\tr(\nu(\alpha+\beta z)))\nonumber\\
& = & \sum_{\nu\in \beta N^*\cap{\cal I}_\pm}e\left(\tr\left(\nu\left(\frac{\alpha}{\beta}+z\right)\right)\right)-
\sum_{{\tiny\begin{array}{c}\nu\in \beta N^*\nonumber\\
\beta'\mu>-\beta\mu'>0\end{array}}}e\left(\tr\left(\nu\left(\frac{\alpha}{\beta}+z\right)\right)\right)\nonumber\\
& = & \sum_{{\tiny\begin{array}{c}\nu\in \beta N^*\\
0<-\nu'<\nu<-(\beta/\beta')\nu'\end{array}}}e\left(\tr\left(\nu\left(\frac{\alpha}{\beta}+z\right)\right)\right).
\label{eq:passaggio}
\end{eqnarray}

The equality (\ref{eq:passaggio}) implies that for every
$\nu\in\spt_K(R_{\alpha,\beta,N})$, \begin{equation}-\frac{\beta}{\beta'}\nu'>\nu>0.\label{eq:etabeta}\end{equation}

Let $\nu$ be as in (\ref{eq:nu_delta}). The condition (\ref{eq:etabeta}) in the 
sum of (\ref{eq:passaggio})
is equivalent to 
\[0<h<\rho(\beta)l,\] because $0<-\nu'<\nu$ holds if and only if (\ref{eq:tonto}) holds, and $\theta>\rho(\beta)$
(thanks to (\ref{eq:111})).

Since $\rho(\beta)\in\QQ_{>0}$, lemma \ref{lemme:rational_functions} applies with $\rho_2=\rho(\beta)$ and 
$\rho_1=0$, and we see
that $R_{0,\beta,N}$ is a rational function
in $\QQ(u,v)$; it is locally analytic because 
$\rho(\beta)>0$. Looking at (\ref{eq:passaggio}), for $\alpha\in K$, we see that:
\[R_{\alpha,\beta,N}(\und{u})=R_{0,\beta,N}\left(\Phi_0\left(\frac{\alpha}{\beta},\frac{\alpha'}{\beta'}\right)\und{u}\right).\]
But the coordinates of $\Phi_0(\alpha/\beta,\alpha'/\beta')\in\TT$ are roots of unit.
Thus, $R_{\alpha,\beta,N}(\und{u})\in\bar{\QQ}(u,v)$ (the coefficients of the Taylor expansion at $\und{0}\in\CC^2$
lie in some
cyclotomic number  field of finite degree over $\QQ$).

The other cases allow a similar proof and are left to the reader. More precisely, there are
three other cases: one has to consider $R_{\alpha,\beta,N}^+$ when $\beta>\beta'>0$, $R_{\alpha,\beta,N}$ when
$\beta>-\beta'>0$ and $R_{\alpha,\beta,N}^+$ when $\beta>-\beta'>0$; in some of these cases,
one has to apply (\ref{eq:222}) too. The
proof of proposition 
\ref{lemme:theta_etc} is complete.

\medskip

\noindent {\bf Remark.} Clearly,
${\cal D}\cap\TT=\Phi_0({\cal W})$ and ${\cal D}^+\cap\TT=\Phi_0({\cal W}^+)$.

\subsubsection{Examples of series $R_{\alpha,\beta,N},R_{\alpha,\beta,N}^+$.}

To let the reader become more familiar with the rational functions introdu\-ced so far,
it is worth to give some explicit examples, even if the proof of theorem 2 can proceed without them.

First of all, if $\beta\in\NN\subset S_+$ or if $\beta\in S\cap\Delta\NN\subset S_\pm$, then, for all $\alpha\in K$:
\[R_{\alpha,\beta,N}(\und{u})=R_{\alpha,\beta,N}^+(\und{u})=0.\] The most important example
is when $\beta$ is an irrational unit $\eta$ of $S$.

Let $\eta>1$ be a unit of $S_+$. By definition we have:
\begin{eqnarray}
f(\und{u}^\eta) & = & f(\und{u})-R_\eta(\und{u})\mbox{ over ${\cal D}$},\label{eq:functional_equation}\\
f^+(\und{u}^\eta) & = & f^+(\und{u})-R_\eta^+(\und{u})\mbox{ over ${\cal D}^+$}\label{eq:functional_equation2}.
\end{eqnarray} 
The identity (\ref{eq:functional_equation}) is the functional equation of $f$, that is
the relation (\ref{eq:fu_eq}) of the introduction. The identity (\ref{eq:functional_equation2}) is a
variant of it, involving the twin series $f^+$. 

Let us inspect more closely these relations.
Let $W=\eta^\ZZ$ be the rank one multiplicative group generated by $\eta$.
We choose the fundamental domain for the multiplicative action of $W$ on 
$K_\pm$
\[{\cal D}(\eta)=\{\nu\in K_\pm:1< -\nu/\nu'\leq\eta^2\},\]
and the fundamental domain for the action of $W$ on $K_+$:
\[{\cal D}_+(\eta)=\{\nu\in K_{+}:\eta^2>\nu/\nu'\geq 1\}.\]
By using the techniques of the proof of {\bf (5)} in proposition \ref{lemme:theta_etc},
it is easily seen that, if
$\und{u}=\Phi_0(z,z')$:
\begin{eqnarray*}
R_\eta(\und{u}) & = & \sum_{\nu\in {\cal D}(\eta)\cap M^*_\pm}e(\tr(\nu z))\\
& = & \sum_{l=1}^{\infty}\sum_{h=1}^{[\rho(\eta)l]}u^lv^h,
\end{eqnarray*}
and
\begin{eqnarray*}
R_\eta^+(\und{u}) & = & \sum_{\nu\in {\cal D}_+(\eta)\cap M^*_\pm}e(\tr(\nu z)),\\
& =& \sum_{l=1}^{\infty}\sum_{h=[\rho^+(\eta)l]+1}^{[\rho^+(1)l]}u^lv^h.
\end{eqnarray*}

By definition of fundamental domain, if $\nu\in K_+$ (resp. $\nu\in K_{\pm}$), then there exists one
and only one  element $k\in\ZZ$ such that $\eta^k\nu\in{\cal D}_+(\eta)$ (resp. $\eta^k\nu\in{\cal D}(\eta)$).

We see that if $k>0$, then
${\cal D}(\eta^k)$ is equal to the disjoint union $\cup_{l=0}^{k-1}\eta^l{\cal D}(\eta)$. This implies:
\[
R_{\eta^k}(\und{u})=\sum_{l=0}^{k-1}R_\eta(\und{u}^{\eta^l}).
\]
In particular, since ${\cal I}_\pm=\lim_{n\rightarrow\infty}{\cal D}(\eta^n)$:
\begin{equation}
f(\und{u})=\sum_{n=0}^\infty R_\eta(\und{u}^{\eta^n})=
\lim_{k\rightarrow\infty}R_{\eta^k}(\und{u}).\label{eq:relation_inf}
\end{equation}
In a similar way, ${\cal D}_+(\eta^k)$ is equal to the disjoint union 
$\cup_{l=0}^{k-1}\eta^l{\cal D}_+(\eta)$ and ${\cal I}_+=\lim_{n\rightarrow\infty}{\cal D}_+(\eta^n)$,
so that:
\begin{equation}
f^+(\und{u})=\sum_{n=0}^\infty R_\eta^+(\und{u}^{\eta^n})=
\lim_{k\rightarrow\infty}R_{\eta^k}^+(\und{u}).\label{eq:relation_inf2}\end{equation}

\subsection{An irrationality criterion for functions.}

We will give equivalent formulations of the property of semi-freeness, by using the linear
independence properties of the rational functions $R_{\alpha,\beta,N},R_{\alpha,\beta,N}^+$, and we will 
prove one half of theorem 2 (the easiest implication). We begin with a useful irrationality criterion.

\begin{Definition}{\em 
A {\em strictly convex cone} of
$\RR^{n}$ of {\em axis} $\und{Y}\in\RR^n$ is the union of all the half-lines $\und{Z}\RR_{>0}$ such that the plane 
angle
between $\und{Y}$ and $\und{Z}$ has absolute value less than, or equal to $\delta$, for some $\delta<\pi/2$.}
\end{Definition}

For a
$2n$-tuple of rational integers \[(p_1,q_1,\ldots,p_n,q_n),\] we consider:
\begin{equation}\gamma_i=\Delta^{-1/2}(p_i-q_i\theta'{}^{-1})\in M^*-\{0\}.
\label{eq:gamma_i_pq}\end{equation} 

Let $Q(\und{V})$ be any formal series in $2n$ variables $\und{V}=(\und{v}_1,\ldots,\und{v}_n)$
with complex coefficients ($n\geq 1$).
We may write:
\[Q(\und{V})=\sum_{\und{p}\in\ZZ^{2n}}
c_{\und{p}}u_1^{p_1}v_1^{q_1}\cdots u_n^{p_n}v_n^{q_n},\]
where $\und{p}=(p_1,q_1,\ldots,p_n,q_n)$ and $c_{\und{p}}\in\CC$. The following definition 
extends definition \ref{def:def_K_spt}.

\begin{Definition}{\em 
The
{\em $K$-support $\spt_K(Q)$ of} $Q$ is the subset of $(\Sigma(M^*))^n\subset\RR^{2n}$ whose 
elements are the $\Sigma^{\oplus n}(\und{\gamma})=(\Sigma(\gamma_1),\ldots,\Sigma(\gamma_n))$'s such that there exists
$\und{p}$ with $c_{\und{p}}\not=0$, satisfying (\ref{eq:gamma_i_pq}).}
\end{Definition}
If $\und{\nu}=(\nu_1,\ldots,\nu_n)\in K^n$ then we will also write $\pi_2(\und{\nu})=(\nu_1',\ldots,\nu_n')\in\RR^n$.
We denote:\[{\goth H}_n^+=K_+\times(K_+\cup K_\pm)^{n-1},\quad{\goth H}_n^\pm=K_\pm\times(K_+\cup K_\pm)^{n-1}.\]
The space of formal series such as $Q(\und{V})$, whose $K$-support is contained in ${\goth H}_n^+$ (resp. 
${\goth H}_n^\pm$), has a structure of $\CC[\und{V}]$-module
(one can multiply a series by a polynomial, but a product on these series is not well defined).
  
\begin{Lemme}[A criterion of irrationality.]
Let $Q(\und{V})$ be any formal series in $2n$ variables 
with complex coefficients, whose $K$-support is contained in the intersection 
of $(M^*)^n$ with ${\goth H}_n^+$ (resp. ${\goth H}_n^\pm$). Let us suppose that 
the set $\pi_2(\spt_K(Q))$ is contained in a strictly convex cone $\Pi$ of $\RR^n$ and that it contains
a sequence of points $(x_{1,s}',\ldots,x_{n,s}')\in(K\setminus\{0\})^n$ such that:
\begin{equation}\lim_{s\rightarrow\infty}(x_{1,s}',\ldots,x_{n,s}')=\und{0}\in\RR^n.\label{eq:y1yn}
\end{equation}
Then $Q$ does not belong to $\CC(\und{V})$.\label{lemme:principle}
\end{Lemme}
\noindent {\bf Proof.} 
We only deal with the case of ${\goth H}_n^+$; the other case allows a very
similar proof, and is left to the reader. 
We choose the strictly convex cone $\Pi$ of the lemma in $\RR_{>0}\times\RR^{n-1}$. 

In order to prove the lemma we need to show that, under our hypotheses, for any non-zero polynomial 
$B\in\CC[\und{V}]$, the $K$-support of $BQ$ is infinite.

Let $B$ be any non-zero polynomial of $\CC[\und{V}]$.
We may write, for complex variables $\und{v}_i=\Phi_0(z_i,z_i')$ with $\und{z}=(z_1,\ldots,z_n)\in\CC^n$,
$\und{z}'=(z_1',\ldots,z_n')\in\CC^n$:
\[B(\und{v}_1,\ldots,\und{v}_n)=\sum_{\und{\lambda}\in{\cal E}}c_{\und{\lambda}}
e(\tr(\und{\lambda}\cdot{}^{{\rm t}}\und{z})),\] where ${\cal E}=\spt_K(B)$
is a non-empty finite subset of $(M^*)^n$, so that the $c_{\und{\lambda}}$'s are all non-zero complex numbers. 

From (\ref{eq:y1yn}) we see that for all $\epsilon>0$, the ball 
$B(\und{0},\epsilon)$ of center $\und{0}\in\RR^n$ and radius
$\epsilon$ contains infinitely many elements $(x_{1,s}',\ldots,x_{n,s}')$ of $\pi_2(\spt_K(Q))$: we denote $B_\epsilon$ the
infinite set of these elements. 

We now prove that $\pi_2(\spt_K(BQ))$ contains a translate (addition of $\RR^n$) 
of $B_\epsilon$ by an element
of $\pi_2({\cal E})$, for $\epsilon>0$ small enough.

The set:
\[Z=\bigcup_{\und{\lambda}\in\pi_2({\cal E})}\und{\lambda}+\pi_2(\spt_K(Q))\]
clearly contains $\pi_2(\spt_K(BQ))$ by the distributive property of product with respect to addition.

Since for all $(y_1,y_2,\ldots,y_n)\in\pi_2(\spt_K(Q))$ we have $y_1>0$ by hypothesis, the finiteness 
of ${\cal E}$ implies that there exists:
\[\iota:=\inf\{y_1\mbox{ such that }(y_1,y_2,\ldots,y_n)\in Z\}.\]
Since $\pi_2(\spt_K(Q))$ has $\und{0}$ as an euclidean adherence point, we have that $\iota\in K$
and there exists a non-empty subset ${\cal E}_0\subset{\cal E}$ such that for all $\und{\lambda}\in\pi_2({\cal E}_0)$,
\[\und{\lambda}=(\iota,*,\ldots,*).\]

Let us choose an element $\und{\lambda}_0\in\pi_2({\cal E}_0)$.
There exists $\epsilon'>0$ such that for all 
$\und{\lambda}\in\pi_2({\cal E})\setminus\{\und{\lambda}_0\}$,
the distance from $\und{\lambda}$ to $\und{\lambda}_0$ is $>\epsilon'$, because $\pi_2({\cal E})$ is finite. 

Thus, the strict convexity of $\Pi$ implies that
for all $\und{\lambda}\in\pi_2({\cal E})\setminus\{\und{\lambda}_0\}$, 
$\und{\lambda}_0\in\RR^n$ is not an adherence point of
$\pi_2(\und{\lambda})+\Pi$, by using the minimality of the first 
coordinate $\iota$ of $\und{\lambda}_0$ and the finiteness
of ${\cal E}$.

Hence, if $0<\epsilon''<\epsilon'$ is small enough, then:
\begin{equation}B(\und{\lambda}_0,\epsilon'')\cap\und{\lambda}+\Pi=\emptyset,\quad\mbox{for all }
\und{\lambda}\in\pi_2({\cal E})\setminus\{\und{\lambda}_0\}.\label{eq:ball_epsilon_2}\end{equation}

On the other side, we already know that the set $B(\und{\lambda}_0,\epsilon'')$
contains the set $\und{\lambda}_0'+B_{\epsilon''}$. Collecting all the informations together,
we have proved, using (\ref{eq:ball_epsilon_2}), that for all $\epsilon''>0$ small enough:
\begin{eqnarray*}
\und{\lambda}_0+B_{\epsilon''} & \subset & \und{\lambda}_0+\pi_2(\spt_K(Q)),\mbox{ and }\\
\und{\lambda}_0+B_{\epsilon''} & \cap & \und{\lambda}+\pi_2(\spt_K(Q))=\emptyset,\mbox{ for all }
\und{\lambda}\in\pi_2({\cal E})\setminus\{\und{\lambda}_0\}.
\end{eqnarray*}
This implies:
\[\und{\lambda}_0+B_{\epsilon''}\subset\pi_2(\spt_K(BQ)).\]
Since $B_{\epsilon''}$ is infinite, $\und{\lambda}_0+B_{\epsilon''}$ is infinite, and the $K$-support of $BQ$
is infinite:
the proof of lemma \ref{lemme:principle} is complete.

\subsection{Linear independence.\label{section:linear}}

Under certain hypotheses, lemma \ref{lemme:principle} becomes an equivalence, and
becomes a tool to check linear independence of functions. In this section we investigate
these properties, and we will prove three results. Lemma
\ref{lemme:rank_one} is a converse of lemma \ref{lemme:principle} under certain hypotheses. 
Lemma \ref{lemme:three_conditions}
makes a connection between linear independence over $\QQ$ of certain functions, and
semi-freeness. Finally, lemma \ref{proposition:equivalence_conditions} resumes the main properties and
is written to apply in one implication of theorem 2.

\medskip

We need some notations.
Let $\alpha_1,\ldots,\alpha_m$ be elements of $K$ and 
$\beta_1,\ldots,\beta_m$ be elements of $S-\{0\}$, such that $\beta_i>0$ for all $i=1,\ldots,m$,
let us denote
$N_i=\beta_i^{-1}M$.

Let us assume, without loss of generality, that there exists $m_0$, with 
$0\leq m_0\leq m$, such that:
\begin{eqnarray*}
\beta_i & \in & S_+\mbox{ for }i=1,\ldots,m_0,\\
\beta_i & \in & S_\pm\mbox{ for }i=m_0+1,\ldots,m.
\end{eqnarray*}
The existence of $m_0$ is guaranteed up to reorder the indexes $i$ (if $m_0=0$ or $m_0=m$, then one of the conditions
is empty).
We write
\begin{eqnarray*}
\Upsilon_i(\und{u}) & = & f(\und{\zeta}_i\und{u}^{\nu_i}),\quad\mbox{ for }i=1,\ldots,m_0,\\
\Upsilon_i(\und{u}) & = & f^+(\und{\zeta}_i\und{u}^{\nu_i})\quad\mbox{ for }i=m_0+1,\ldots,m,
\end{eqnarray*} so that
for all $i=1,\ldots,m$, we have \[\spt_K(\Upsilon_i)\subset N_i^*\cap K_\pm,\quad i=1,\ldots,m.\]
Together with the functions $\Upsilon_i$, we need to also consider:
\begin{eqnarray*}
\Upsilon_i^+(\und{u})&=&f^+(\und{\zeta}_i\und{u}^{\nu_i}),\quad\mbox{ for }i=1,\ldots,m_0,\\
\Upsilon_i^+(\und{u})&=&f(\und{\zeta}_i\und{u}^{\nu_i}),\quad\mbox{ for }i=m_0+1,\ldots,m,
\end{eqnarray*} so that 
\begin{equation}
\spt_K(\Upsilon_i^+)\subset N_i^*\cap{\cal I}_+,\quad i=1,\ldots,m.\label{eq:Upsilon_plus}
\end{equation}
Let us also define:
\begin{eqnarray}
Q(\und{u})&=&\sum_{i=1}^mc_i\Upsilon_i(\und{u}),\label{eq:Qpm}\\
Q^+(\und{u})&=&\sum_{i=1}^mc_i\Upsilon_i^+(\und{u}).\nonumber
\end{eqnarray} 

\begin{Lemme}
Let $\beta\in S-\{0\}$. The map $\TT\rightarrow\TT$ defined by $\und{u}\mapsto\und{u}^\beta$
is a group homomorphism whose kernel is a finite subgroup of $\TT$ with $|\no(\beta)|=|\beta\beta'|$
elements.

Let $\und{\zeta}\in\TT$ be a torsion point. There exists an irrational unit $\eta\in S_+$
such that $\und{\zeta}^\eta=\und{\zeta}$.
\label{lemma:torsion_3}\end{Lemme}
\noindent {\bf Proof.} 
The first part of the lemma is clear, since $\beta^{-1}M/M\cong M^*/(\beta M^*)$ is a 
finite group with $|\no(\beta)|$ elements. We call such a kind of map an {\em isogeny of degree} $|\no(\beta)|$.

Let $\eta$ be a unit of $S$: the isogeny $\TT\rightarrow\TT$ defined by $\und{u}\mapsto\und{u}^\eta$
is an {\em automorphism} (its degree is $1$).

The set $\und{\zeta}^S=\{\und{\zeta}^\beta;\beta\in S\}$
is a finite subgroup of $\TT$ and is stable under the action of $S$. The automorphisms constructed
with units of $S$ act as permutations of the finite set $\und{\zeta}^S$, so that
a suitable non-zero integer power of any given irrational unit of $S$ induces the identity map on $\und{\zeta}^S$.

\subsubsection{A criterion of linear independence over $\CC$.}

Let $\eta>1$ be a unit of $S_+$ such that
$\Phi_0(\alpha_i\eta,\alpha_i'\eta')=\Phi_0(\alpha,\alpha')$ for all $i=1,\ldots,m$ (whose existence is 
guaranteed by lemma
\ref{lemma:torsion_3}).

\begin{Lemme} Let $c_1,\ldots,c_m$ be complex numbers. The following conditions are equivalent:
\begin{itemize}
\item[(1)] The function $Q^+(\und{u})$ in (\ref{eq:Qpm}) is
rational.
\item[(2)] There exists a positive real number $\epsilon$, depending only on $Q^+$, such 
that for all $\alpha\in\pi_2(\spt_K(Q^+))$,
we have that
$\alpha>\epsilon$.
\item[(3)] We have the identity of functions: \begin{equation}
\sum_{i=1}^{m}c_iR_{\eta,N_i}^+(\und{\zeta}_i^\sharp\und{u})=0,\label{eq:lin_dep_rat}\end{equation}
where $N_i=\beta_i^{-1}M$ and $\und{\zeta}_i^\sharp=\Phi_0(\alpha_i/\beta_i,\alpha_i'/\beta_i')$.
\item[(4)] The function $Q(\und{u})$ in (\ref{eq:Qpm}) is rational.
\item[(5)] There exists a positive real number $\omega$, depending only on $Q$, such 
that for all $\alpha\in\pi_2(\spt_K(Q))$,
we have that
$\alpha<-\omega$.
\item[(6)] We have the identity of functions: \begin{equation}
\sum_{i=1}^{m}c_iR_{\eta,N_i}(\und{\zeta}_i^\sharp\und{u})=0,
\label{eq:lin_dep_rat2}
\end{equation}
where $N_i$ and $\und{\zeta}_i^\sharp$ are defined as in the point (3).
\end{itemize}
\label{lemme:rank_one}
\end{Lemme}
{\bf Proof.} We first prove that (3) implies (1).
From (\ref{eq:functional_equation2}), or 
(\ref{eq:relation_inf}), we see that:
\[\sum_{k=0}^\infty R_{\eta,N_i}^+(\und{\zeta}_i^\sharp
\und{u}^{\eta^k})=f^+_{N_i}(\und{\zeta}_i^\sharp\und{u}),\quad i=1,\ldots,m,\]
because $W=\eta^\ZZ$ fixes the torsion
points $\und{\zeta}_i$. 
Hence, if (\ref{eq:lin_dep_rat}) holds, then
\[\sum_{i=1}^{m}c_i\sum_{k=0}^\infty R_{\eta,N_i}^+(\und{\zeta}_i^\sharp\und{u}^{\eta^k})=
\sum_{i=1}^mc_if^+_{N_i}(\und{\zeta}_i^\sharp\und{u})=0.\] 
We have, for all $i=1,\ldots,m$, that $f_{N_i}^+(\und{\zeta}_i^\sharp\und{u})-\Upsilon_i(\und{u})=
R_{\alpha_i,\beta_i}^+(\und{u})$ is rational because of proposition \ref{lemme:theta_etc}.
Thus \[\sum_{i=1}^mc_i\Upsilon_i^+(\und{u})=-\sum_{i=1}^mc_iR_{\alpha_i,\beta_i}^+(\und{u})\]
is rational.

\vspace{10pt}

That (1) implies (2) follows by the lemma \ref{lemme:principle} with $n=1$;
indeed, the $K$-support of $Q^+(\und{u})$ is contained in ${\goth H}_n^+=M^*_+$ (in this case, the 
strictly convex cone is a half-line).
As $Q^+(\und{u})$ is rational, $0$ is not an 
adherence point of $\pi_2(\spt_K(Q^+))$ and the required positive real number $\epsilon$ exists. One can also apply
lemma \ref{lemme:rational_functions}.
 
\vspace{10pt}

We now show that (2) implies (3). For every $\gamma\in\spt_K(Q^+)$, \[\gamma'>\epsilon\] for
some $\epsilon>0$, by hypothesis.
From the definition of the functions $R_{\alpha,\beta},R_{\alpha,\beta}^+$, we obtain the 
equality:
\begin{equation}
\sum_{i=1}^mc_if^+_{N_i}(\und{\zeta}_i^\sharp\und{u})=Q^+(\und{u})+R(\und{u}),
\label{eq:supports_bornes}\end{equation}
where \[R(\und{u})=\sum_{i=1}^{m}c_iR^+_{\alpha_i,\beta_i}(\und{u})\] is a rational function,
by proposition \ref{lemme:theta_etc}. The point {\bf (4)} of proposition \ref{lemme:theta_etc}
says, with lemma \ref{lemme:theta_etc} (we take ${\goth H}_1^\pm=K_\pm$), that there exists $\epsilon'>0$ such that for all 
$\gamma\in\pi_2(\spt_K(R))$, \[\gamma>\epsilon'.\]

Let ${\cal Q}$ be the $K$-support of $Q^++R$. We obtain from above, that for all $\gamma\in\pi_2({\cal Q})$,
\[\gamma>\epsilon'',\] with
$\epsilon''=\min(\epsilon,\epsilon')$.

Let ${\cal P}$ be the $K$-support of 
the rational function
\[P(\und{u})=\sum_{i=1}^{m}c_i R_{\eta,N_i}^+(\und{\zeta}_i^\sharp\und{u}).\] 
We want to prove that ${\cal P}$ is empty; this will imply (\ref{eq:lin_dep_rat}).

From (\ref{eq:supports_bornes}), we see that: 
\[{\cal Q}=\bigcup_{i=0}^\infty\eta^i{\cal P}\] (disjoint union).
Suppose by contradiction that there exists some $\nu\in{\cal P}$.
Then for all $i\geq 0$ we have that $\eta^i\nu\in{\cal Q}$. As \[\lim_{i\rightarrow\infty}(\eta^i\nu)'=0,\]
for some $i$ big enough we have that \[\epsilon''>(\eta^i\nu)'>0,\] and this implies that there
exists some element $\gamma\in{\cal Q}$ such that $0<\gamma'<\epsilon''$. But $\gamma$ does not
belong to the $K$-support of $R(\und{u})$ because of our choice of $\epsilon''$, and this
means that $\gamma\in\spt_K(Q^+(\und{u}))$: this gives us the required contradiction.
Assuming part (2) of the lemma, it is now clear that (\ref{eq:lin_dep_rat})
holds.

\vspace{10pt}

Let us now prove that (1) implies (4). The point {\bf (3)} of proposition
\ref{lemme:theta_etc} implies that \[\Upsilon_i(\und{u})=
-\Upsilon_i^+(\und{u})+R_i(\und{u})\] for some  rational functions
$R_i(\und{u})\in \bar{\QQ}(\und{u})$. Thus:
\[Q(\und{u})=-Q^+(\und{u})+\sum_{i=1}^mc_iR_i(\und{u})\] is a rational function in $\CC(\und{u})$.
The proof that (4) implies (1) is similar.

\vspace{10pt}

The proof that (4), (5), (6) are equivalent, runs along the same ideas than the proof for the equivalence of (1), (2), (3).
It is enough to remark that $Q(\und{u})$ has its $K$-support contained in $M^*_\pm$. 
Then one applies the lemma \ref{lemme:principle} with $n=1$ and ${\goth H}_n^\pm=M^*_\pm$:
the proof of the lemma \ref{lemme:rank_one} is complete.

\subsubsection{Linear independence over $\QQ$.}

\begin{Definition}
{\em Let $N\supset M$ be a complete $\ZZ$-module, let $\und{u}$ and $(z,z')$ be couples of
variables such that $\und{u}=\Phi_0(z,z')$. 
The {\em Hecke's
geometric series}
$A_N,B_N$ are defined by:
\begin{eqnarray*}
A_N(\und{u})&=&\sum_{\nu\in N_\pm^*}e(\tr(\nu z)),\\
B_N(\und{u})&=&\sum_{\nu\in N_+^*}e(\tr(\nu z))
\end{eqnarray*}}
\label{def:HGS}\end{Definition}
(see \cite{Hecke:Analytische}). The series $A_N$ converges over ${\cal H}\times{\cal H}^-$, where
${\cal H}^-$ is the lower half-plane
\[{\cal H}^-=\{z\in\CC\mbox{ such that }\Im(z)<0\},\]
and the series $B_N$ converges over ${\cal H}\times{\cal H}$.
\begin{Lemme}
Let $c_1,\ldots,c_m$ be rational numbers. The following conditions are equivalent.
\begin{enumerate}
\item The linear relation 
\begin{equation}\sum_{i=1}^mc_i\FF_{N_i}\left(\Phi_0\left(\frac{\alpha_i}{\beta_i},\frac{\alpha_i'}{\beta_i'}
\right)\und{U}\right)=0\label{relation}\end{equation}
holds.
\item We have the linear dependence relations
of Hecke geometric series:
\begin{eqnarray}
\sum_{i=1}^mc_iA_{N_i}\left(\Phi_0\left(\frac{\alpha_i}{\beta_i},\frac{\alpha_i'}{\beta_i'}
\right)\und{u}\right) & =
&0\label{eq:relationH0} \\ 
\sum_{i=1}^mc_iB_{N_i}\left(\Phi_0\left(\frac{\alpha_i}{\beta_i},\frac{\alpha_i'}{\beta_i'}\right)\und{u}\right) &=
&0.\label{eq:relationH1}\end{eqnarray}
\item We have the linear dependence relations
of rational functions:
\begin{eqnarray}
\sum_{i=1}^mc_iR_{\eta,N_i}\left(\Phi_0\left(\frac{\alpha_i}{\beta_i},\frac{\alpha_i'}{\beta_i'}
\right)\und{u}\right) & =
&0\label{eq:relation0} \\ 
\sum_{i=1}^mc_iR_{\eta,N_i}^+\left(\Phi_0\left(\frac{\alpha_i}{\beta_i},\frac{\alpha_i'}{\beta_i'}\right)\und{u}\right) &=
&0.\label{eq:relation1}\end{eqnarray}
\end{enumerate}
\label{lemme:three_conditions}
\end{Lemme}

\noindent {\bf Proof.} If $N$ is a complete $\ZZ$-module of $K$, we have that: \[N^*=N^*_\pm\cup
N^*_+\cup(-N^*_\pm)\cup(-N^*_+)\] (disjoint union). We have the identity of formal series:
\begin{equation}
\FF_N(\und{U})=A_N(\und{U})+B_N(\und{U})+C_N(\und{U})+D_N(\und{U}),\label{eq:identity_formal}\end{equation}
where $A_N,B_N$ are the series of
definition \ref{def:HGS}, and
\begin{eqnarray*}
C_N(\und{u})&=&\sum_{\nu\in -N_\pm^*}e(\tr(\nu z)),\\
D_N(\und{u})&=&\sum_{\nu\in -N_+^*}e(\tr(\nu z))
\end{eqnarray*}
($C_N$ converges over ${\cal H}^-\times{\cal H}$, and $D_N$ converges over ${\cal H}^-\times{\cal H}^-$).
Thus, part 1 of the lemma implies part 2.

Moreover, if $N=\beta^{-1}M$ for $\beta\in S\setminus\{0\}$, 
we have that $N^*_\pm=\cup_{k\in\ZZ}(\eta^k{\cal D}(\eta)\cap N^*)$ (disjoint
union), and similarly for 
$N^*_+$. Thus:
\begin{eqnarray*}
\sum_{k\in\ZZ}R_{\eta,N}(\und{u}^{\eta^k})&=&A_N(\und{u}),\\
\sum_{k\in\ZZ}R_{\eta,N}^+(\und{u}^{\eta^k})&=&B_N(\und{u}).
\end{eqnarray*}
If $\alpha\in K,\beta\in S\setminus\{0\}$, and $\eta$ is a
unit such that $\alpha\eta-\alpha\in M$, then for all
$k\in\ZZ,\nu\in \beta M^*$:
\[\tr(\nu\alpha\eta^k/\beta)-\tr(\nu\alpha/\beta)\in\ZZ.\]
In particular, for all $i=1,\ldots,m$, the hypotheses of the lemma imply:
\begin{equation}
e(\tr(\nu\alpha_i\eta^k/\beta_i))=e(\tr(\nu\alpha_i/\beta_i)).
\label{eq:nu_eta}\end{equation}
Now, by using (\ref{eq:nu_eta}) we have, for $\und{u}=\Phi_0(z,z')$:
\begin{eqnarray*}
A_{N_i}\left(\Phi_0\left(\frac{\alpha_i}{\beta_i},\frac{\alpha_i'}{\beta_i'}\right)
\und{u}\right) & = &
\sum_{\nu\in (N_i^*){}_\pm}e(\tr(\nu\alpha_i/\beta_i))e( \tr(\nu z))\\
&=&\sum_{k\in\ZZ}R_{\eta,N}\left(\Phi_0\left(\frac{\alpha_i}{\beta_i},\frac{\alpha_i'}{\beta_i'}\right)^{\eta^k}
\und{u}^{\eta^k}\right),
\end{eqnarray*}
and similarly for $B_{N_i}$.
From these identities, it follows that part 2 of the lemma is equivalent to part 3.

We prove that part 2 of the lemma implies part 1.
We have: \begin{eqnarray*}
\overline{A_{N_i}\left(\Phi_0\left(\frac{\alpha_i}{\beta_i},\frac{\alpha_i'}{\beta_i'}\right)
\und{u}\right)} & = &
\sum_{\nu\in (N_i^*){}_\pm}\overline{e(\tr(\nu\alpha_i/\beta_i))e( \tr(\nu z))}\\
& = & \sum_{\nu\in (N_i^*){}_\pm}e(\tr(-\nu\alpha_i/\beta_i))\overline{e( \tr(\nu z))}\\
& = & \sum_{\nu\in (N_i^*){}_\pm}e(\tr(-\nu\alpha_i/\beta_i))e( \tr(-\nu \overline{z}))\\
& = & \sum_{\mu\in (-N_i^*){}_\pm}e(\tr(\mu\alpha_i/\beta_i))e( \tr(\mu \overline{z}))\\
& = &
C_{N_i}\left(\Phi_0\left(\frac{\alpha_i}{\beta_i},\frac{\alpha_i'}{\beta_i'}\right)\und{u}^*\right),
\end{eqnarray*}
with $\und{u}^*=\Phi_0(\overline{z},\overline{z}')$ (the symbol $\overline{\cdot}$ means Òcomplex conjuguation"). Similarly,
one obtains:
\begin{eqnarray*}
\overline{B_{N_i}\left(\Phi_0\left(\frac{\alpha_i}{\beta_i},\frac{\alpha_i'}{\beta_i'}\right)
\und{u}\right)} & = & D_{N_i}\left(\Phi_0\left(\frac{\alpha_i}{\beta_i},\frac{\alpha_i'}{\beta_i'}\right)\und{u}^*\right).
\end{eqnarray*}
Finally, the linear relations 
\begin{equation}
\sum_{i=1}^mc_iX_i\left(\Phi_0\left(\frac{\alpha_i}{\beta_i},\frac{\alpha_i'}{\beta_i'}
\right)\und{U}\right)=0,\label{eq:ABCD}\end{equation}
for $X_i=A_{N_i},B_{N_i}$ also imply the linear relations (\ref{eq:ABCD}) for $X_i=C_{N_i},D_{N_i}$, because the coefficients
$c_i$ are rational integers, thus invariant by complex conjugation: hence the relations
(\ref{eq:ABCD}) hold for $X_i=A_i,B_i,C_i,D_i$, and applying (\ref{eq:identity_formal}),
we obtain the formal relations (\ref{relation}). The proof of lemma \ref{lemme:three_conditions} is complete.

\medskip

\noindent {\bf Remark.} The series $A_N,B_N,C_N,D_N$ converges, but have disjoint domains of convergence.

\medskip

\begin{Lemme}
Let $c_1,\ldots,c_m$ be rational numbers.
Taking into account the notations introduced above, the following conditions are equivalent.
\begin{enumerate}
\item We have that $Q(\und{u})$, as in (\ref{eq:Qpm}), belongs to $\bar{\QQ}(\und{u})$. 
\item The linear dependence relation 
(\ref{relation}) holds.
\end{enumerate}
\label{proposition:equivalence_conditions}\end{Lemme}

\noindent {\bf Proof.}
By lemma \ref{lemme:three_conditions}, the second condition of the lemma \ref{proposition:equivalence_conditions} 
implies (\ref{eq:relation1}). This identity 
is the condition (6) of lemma \ref{lemme:rank_one}) which is equivalent to condition (4) which is
the first condition of lemma \ref{proposition:equivalence_conditions}. 

On the other side, if the first condition of the lemma holds, by lemma 
\ref{lemme:rank_one}, conditions (3) and (6) of lemma
\ref{lemme:rank_one} hold so that (\ref{eq:relation0}) and
(\ref{eq:relation1}) are satisfied; then one applies lemma \ref{lemme:three_conditions} to obtain 
the second condition of the lemma.

\subsection{Proof of one half of theorem 2.\label{section:portrait2}}

First of all, we state and prove a little lemma.

\begin{Lemme}
The set ${\cal D}$ does not contain any torsion point of $\TT$.
\label{lemme:no_torsion}\end{Lemme}
\noindent {\bf Proof}. If $(u,v)=\Phi_0(z,z')\in\TT\cap{\cal D}$, then
\begin{equation}(z,z')\not\in\RR\times\CC\label{eq:condition_u_4},\end{equation}
because $(z,z')\in\RR\times\CC^\times$ if and only if $|u|\not=1$ and $|u||v|^\theta=1$,
and $(z,z')\in\RR\times\{0\}$ if and only if $|u|=|v|=1$, by (\ref{eq:more_explicit}).

In particular by (\ref{eq:condition_u_4}), if $(u,v)\in\TT\cap{\cal D}$ then $z\not\in \Sigma(K)$ and
$\und{u}$ is not a torsion point: every point of ${\cal D}$ is a point of infinite order.

\medskip

We prove here the easiest implication of theorem 2.
\begin{Proposition}
Let ${\cal
M}=(\und{u}_1,\ldots,\und{u}_m)$ be a $m$-tuple of algebraic elements of
$\TT(\bar{\QQ})\cap{\cal D}$ which is not semi-free. then there exists
$c_1,\ldots,c_m$ rational numbers not all zero, and an algebraic number
$\lambda$ such that \[\sum_{i=1}^mc_if(\und{u}_i)=\lambda.\]
\label{proposition:one_half}\end{Proposition} 

\noindent {\bf Proof.} By lemma \ref{lemme:no_torsion} we know that $\und{u}_i$ has infinite order for all $i$.
Without loss of generality, we can suppose that there exist elements
$\alpha_1,\ldots,\alpha_m$ of $K$, elements 
$\beta_1,\ldots,\beta_m$ of $S-\{0\}$ such that 
\begin{equation}\und{u}_i=\Phi_0(\alpha_i,\alpha_i')\und{v}^{\beta_i}\label{eq:torsion_and_not_2}\end{equation}
for an element $\und{v}\in\TT(\bar{\QQ})$ of infinite order. 

In view of the semi-freeness property, we can proceed for the most comfortable possible choice of $\alpha_i$ and of the point $\und{v}$
in (\ref{eq:torsion_and_not_2})
(by using the lemmata \ref{lemma:not_depend} and \ref{lemme:two_groups} of the appendix).

In particular, we may modify (\ref{eq:torsion_and_not_2}) so that $\und{v}$ is replaced with $\und{v}^{\pm \eta^k}$
for some unit $\eta$, consequently, $\beta_i$ is replaced by $\pm\beta_i\eta^{-k}$.

Up to replace $\und{v}$ by $\und{v}^{-1}$, we may suppose that $\und{v}=\Phi_0(z,z')$ with $(z,z')\in{\cal H}\times\CC$
(from lemma \ref{lemme:no_torsion}, or (\ref{eq:condition_u_4}) we know that $(z,z')\not\in\RR\times\CC$). 

Since $\und{v}\in\Phi_0({\cal H}\times\CC)$, then for $k\in\NN$ big enough: \[\und{v}^{\eta^k}\in{\cal D}^+.\] Indeed, 
$\und{v}^{\eta^k}=\Phi_0(\eta^kz,\eta'{}^kz')$; since $\eta'{}^k=\eta^{-k}$, we have that
$(\eta^kz,\eta'{}^kz')\in{\cal W}^+$ for $k$ big enough. This means that we can choose $\und{v}\in{\cal D}^+$
in (\ref{eq:torsion_and_not_2}); from these
equalities we also see that $\beta_i>0$ for all $i=1,\ldots,m$.

Taking into account (\ref{eq:torsion_and_not_2}) and applying lemma \ref{proposition:equivalence_conditions},
we see that (\ref{eq:Qpm}) holds, with rational coefficients $c_i$
not all vanishing, so that $Q(\und{u})\in\bar{\QQ}(\und{u})$. By definition of $\Theta_M$:
\[\Upsilon_i(\und{u})=\Theta_M(\und{\zeta}_i^{\sharp}\und{u}^{\beta_i})-f(\und{\zeta}_i^{\sharp}\und{u}^{\beta_i}),\quad
i=m_0+1,\ldots,m,\] so that  
\[\sum_{i=1}^{m_0}c_if(\und{\zeta}_i^{\sharp}\und{u}^{\beta_i})-\sum_{j=m_0+1}^{m}c_jf(\und{\zeta}_j^{\sharp}\und{u}^{\beta_j})
=Q(\und{u})-\sum_{j=m_0+1}^mc_j\Theta_M(\und{\zeta}_j^{\sharp}\und{u}^{\beta_j})\]
is a rational function in $\bar{\QQ}(\und{u})$ (because $\Theta_M\in\QQ(\und{u})$).
The latter equality holds for $\und{u}=
\und{v}$ because $\und{v}\in{\cal D}^+$, and all the series above converge on ${\cal D}^+$. Hence,
\[\sum_{i=1}^{m_0}c_if(\und{u}_i)-\sum_{j=m_0+1}^{m}c_jf(\und{u}_j)
=\lambda\in\bar{\QQ},\] and the proof of proposition \ref{proposition:one_half} is complete.

\medskip

In the rest of this article we prove the other implication of the theorem 2, that is: if 
$\und{u}_1,\ldots,\und{u}_m\in\TT(\bar{\QQ})\cap{\cal D}$ are such that the numbers
$f(\und{u}_1),\ldots,f(\und{u}_m)$ 
are algebraically dependent over $\QQ$, then these points of $\TT(\bar{\QQ})$ are
multiplicatively semi-independent (with respect to $\Phi_0$). 
The existence of a linear relation such as (\ref{eq:relations_1}) will follow from proposition 
\ref{proposition:one_half}.

\section{Arithmetic and analytic tools.\label{section:tools}}
In this section we deal with the arithmetic and analytic tools to be used in the proof of our theorem.
We start with the (arithmetic) result which is the heart of the Mahler method used here; the
criterion of algebraic independence of Loxton and van der Poorten. 

To apply this criterion, we  have to
check some technical hypotheses, notably a certain analytic 
condition (property A); we first study in more detail the algebraicity of the 
action of $S$ over $\TT$ defined above, as well as the nature of the 
automorphism involved in the functional equation of $f$. Then we relate these algebraic
properties to the Òanalytic properties" required by the Mahler method, by using 
Masser's vanishing theorem.
Particular bases for $S$-sub-modules of $\TT$ are constructed. We do this
by extending a technical lemma (lemma 3 p. 37 of \cite{Loxton:Fredholm}) of Loxton and van der Poorten (\footnote{The
generalisation needs Baker's theorem on linear forms of logarithms, unlike the original lemma of Loxton
and van der Poorten.}). 

Then we present the 
functions in two variables to be used, as well as some of their elementary properties.
Finally, we are ready to introduce the locally analytic functions in several variables
which will be studied to prove our theorem 2, and the functional equation satisfied by these
functions; here, the results of section \ref{section:semifree} will be applied.

\subsection{Arithmetic tools: an outline of Mahler's theory.}

Following \cite{Masser:Hecke}, we recall a terminology that will be used in the sequel.
Let \[{\cal B}=(a_{i,j})_{1\leq i,j\leq \tilde{n}}\] be a regular square matrix of order $\tilde{n}$ 
with its coefficients $a_{i,j}$ in $\ZZ$,
Let $\und{V}=(v_1,\ldots,v_{\tilde{n}})$ be an element of $\GG_m^{\tilde{n}}(\CC)$.
We write \[{\cal B}.\und{V}=\tilde{\und{V}}\in\CC^{\tilde{n}},\] where 
$\tilde{\und{V}}=(\tilde{v}_1,\ldots,\tilde{v}_{\tilde{n}})$ with 
$\tilde{v}_i=\prod_{j=1}^{\tilde{n}}v_j^{a_{i,j}}$. 

\begin{Definition}{\em We call
{\em isogeny associated to ${\cal B}$} the map $\GG_m^{\tilde{n}}(\CC)\rightarrow\GG_m^{\tilde{n}}(\CC)$ 
given by $\und{V}\mapsto{\cal B}.\und{V}$. }\end{Definition}

Suppose that for all $1\leq i,j\leq \tilde{n}$ we have $a_{i,j}\geq 0$.
It is well known (p. 396 of \cite{Loxton:Arithmetic}) that the maximum 
$\lambda_{{\cal B}}$ of the absolute values of the eigenvalues of ${\cal B}$ 
is itself an eigenvalue of ${\cal B}$.

\begin{Definition}{\em We say that the matrix ${\cal B}$ is {\em good} if it is non-singular, it has no roots of
unit as eigenvalues, and
it has an eigenvector $\und{v}_{{\cal B}}\in\RR^{\tilde{n}}\subset\CC^{\tilde{n}}$ corresponding to
$\lambda_{{\cal B}}$ whose coordinates are all positive.
(this is the definition of \cite{Masser:Hecke} p. 209).}\end{Definition}

If ${\cal B}$ is good then there exists a non-empty subset 
${\cal U}({\cal B})$ of $\GG_m^{\tilde{n}}(\CC)$ 
such that for all $\und{U}\in{\cal U}({\cal B})$ we have
\[\lim_{k\rightarrow\infty}{\cal B}^k.\und{U}=\und{0}\]
in $\CC^{\tilde{n}}$. 

The set ${\cal U}({\cal B})$ is moreover the intersection of an open euclidean neighbourhood
of $\und{0}$ in $\CC^{\tilde{n}}$ with
$\GG_m^{\tilde{n}}(\CC)$
(see the definition 2 p. 93 of \cite{Loxton:Variables} or the definition
2 p. 397 of \cite{Loxton:Arithmetic} and 
the lemma 1 p. 397 of \cite{Loxton:Arithmetic}).

\begin{Definition}{\em If $\und{A}\in{\cal U}({\cal B})$, we say that $\und{A}$
{\em satisfies the property A} (abridged expression for Òanalytic property") if the only locally analytic
function $F$ such that 
$F({\cal B}^k.\und{A})=0$ for all $k$ big enough is the zero function (p. 398 of
\cite{Loxton:Arithmetic}).}\end{Definition}

Let $\Psi_1(\und{V}),\ldots,\Psi_m(\und{V})$ be 
locally analytic functions in $\tilde{n}$ variables and satisfying
a system of functional equations:
\begin{equation}
\Psi_j({\cal B}.\und{V})=\Psi_j(\und{V})+R_j(\und{V})\mbox{ for }
1\leq j\leq m,\label{eq:eqfonctionnelles}\end{equation}
where $R_1,\ldots,R_m$ are rational functions.
We will use the theorem on p. 399 of 
\cite{Loxton:Arithmetic} that we quote here as a proposition.

\vspace{10pt}

\begin{Proposition}[Loxton and van der Poorten.]
If the Taylor series at $\und{0}$ of $\Psi_j$, the coefficients of $R_j$ and
the coordinates of $\und{A}$ are algebraic numbers for all $i=1,\ldots,m$, if 
$\und{A}\in{\cal U}({\cal B})$ satisfies the property A,
if $\Psi_i$ and $R_i$ are defined and analytic in a neightbourhood of $\und{A}$, and if
\[\Psi_1(\und{V}),\ldots,\Psi_m(\und{V})\] are
algebraically independent over $\CC(\und{V})$, then the complex numbers
\[\Psi_1(\und{A}),\ldots,\Psi_m(\und{A})\] are algebraically independent over $\QQ$.
\label{propo:loxton}\end{Proposition}

\subsection{Isogenies.\label{section:analytic_tools}}

Let $\nu$ be an element of $S\setminus\{0\}$, let us denote:
\begin{equation}
{\cal B}_0(\nu)={\goth B}_0\cdot\sqm{\nu}{0}{0}{\nu'}\cdot{\goth B}_0^{-1}.\label{eq:explicit}
\end{equation}
As for (\ref{eq:explicit_1}), ${\cal B}_0(\nu)$ has determinant $\no(\nu)$ and 
has rational integer coefficients; moreover, it satisfies $\und{u}^\nu={\cal B}_0(\nu).\und{u}$ for the choice 
of the exponential function $\Phi_0$.

\begin{Lemme}
The isogeny $\und{u}\mapsto\und{u}^\nu$ of $\TT$ extends to an analytic
map $\CC^2\rightarrow\CC^2$ if and only if $\Sigma(\nu)\in{\cal X}$ where:
\[{\cal
X}=\left\{(y,y')\in\RR^2\mbox{
such that }y\geq \max\left\{y',\frac{\theta'}{\theta}y'
\right\}>0\right\}.\]\label{lemma:set_S}\end{Lemme}

\noindent {\bf Proof.} This can be easily checked by a direct computation,
because the first condition is equivalent to require that all the entries in the left 
hand side of (\ref{eq:explicit}) are
$\geq 0$; if $\Sigma(\nu)\in{\cal X}$, then all of the entries of ${\cal B}_0(\nu)$
are non-negative rational integers. 

\medskip

\noindent {\bf Example.} If $\eta\in S_+$ is a unit such that $\eta>1$, then $\Sigma(\eta)\in{\cal X}$,
so that the isogeny (automorphism of infinite order)
\begin{eqnarray*}
\TT & \rightarrow & \TT\\
\und{u}& \mapsto& \und{u}^\eta
\end{eqnarray*}
is locally analytic. 

This can also be seen by using continued fractions as follows. Using (\ref{eq:continued_fraction}) of the appendix we get:
\[{\cal B}_0(\eta)=\sqm{b_{2p-1}}{1}{1}{0}\cdot\cdots\cdot\sqm{b_{0}}{1}{1}{0},\] which implies that
all the entries of ${\cal B}_0(\eta)$ are non-negative integers. Clearly ${\cal B}_0(\eta)$
is good.

\subsubsection{Direct sums of isogenies.}

In the following, we must work in finite direct sums of copies of $\TT$, and we also need to 
consider direct sums of
isogenies. Let \[\Sigma^{\oplus n}:K^n\rightarrow\RR^{2n}\] be the direct sum
of $n$ copies of the embedding $\Sigma$.

\begin{Lemme} Let $\und{\nu}:=(\nu_1,\dots,\nu_n)\in (S\setminus\{0\})^n$. Then the isogeny
$\und{V}\mapsto\und{V}^{\und{\nu}}$ defined by:
\begin{eqnarray*}
\TT^n & \rightarrow & \TT^n\\
(\und{v}_1,\ldots,\und{v}_n) & \mapsto & (\und{v}_1^{\nu_1},\ldots,\und{v}_n^{\nu_n})
\end{eqnarray*}
extends to an analytic map $\TT^n\rightarrow\TT^n$ if and only if
$\Sigma^{\oplus n}(\und{\nu})\in{\cal X}^n.$\label{lemme:action_S}\end{Lemme}
\noindent {\bf Proof.} This is a simple consequence of lemma
\ref{lemma:set_S}.

\medskip

\noindent {\bf Example}. The isogeny
\[\und{V}=(\und{v}_1,\ldots,\und{v}_n)\mapsto
\und{V}^{\und{\eta}}:=(\und{v}_1^\eta,\ldots,\und{v}_n^\eta)\]
extends to an analytic map $\CC^{2n}\rightarrow\CC^{2n}$.
Setting 
\[{\cal B}:={\cal B}_0(\eta)^{\oplus n}\]
we see that ${\cal B}$ is good, and that
\[\und{V}^{\und{\eta}}={\cal B}.\und{V}.\]
\subsection{Multiplicative independence.}
\begin{Definition} {\em We say that the points $\und{a}_1,\ldots,\und{a}_n\in\TT$
are {\em multiplicatively independent (with respect to $\Phi_0$)} 
if the only solution $(\mu_1,\ldots,\mu_n)\in S^n$ of
\begin{equation}\und{a}_1^{\mu_1}\cdots\und{a}_n^{\mu_n}=\und{1}
\label{eq:relation111}\end{equation}
in $\TT$ is the trivial solution $\mu_1=\cdots=\mu_n=0$.
Clearly a notion of multiplicative dependence is determined as well.}\label{definition:mult}\end{Definition}

\noindent {\bf Remark.} It is easy to see, applying directly the definition \ref{defi:semi_free},
that if $\und{a}_1,\ldots,\und{a}_n$
are multiplicatively independent, then they are also multiplicatively semi-independent.
The converse is false, and the section \ref{section:portrait} of the appendix 
provides some counterexamples.

\subsubsection{$S$-groups of finite rank.}

As $S$ may have non-principal ideals, it is in general difficult to study
abstract $S$-modules. The definition below will be useful in the following.

\begin{Definition}{\em Let $n$ be a non-negative integer:
we say that an $S$-module $\Gamma$ is an $S$-{\em group of rank }$n$
if \[\Gamma\cong\Gamma_{{\tiny\mbox{ tors}}}\oplus \underbrace{S\oplus\cdots\oplus S}_{\mbox{ $n$ times}}\]
where $\Gamma_{{\tiny\mbox{ tors}}}$ is a finite group (in case $n=0$, the definition requires
that $\Gamma=\Gamma_{{\tiny\mbox{ tors}}}$).}\end{Definition}

The following lemma provides a tool to simplify matters in case we only consider
finitely generated $S$-submodules of $\TT$ (this will be the case all along this article).
\begin{Lemme}
Any finitely generated $S$-module $\Lambda\subset\TT$ is contained in some $S$-group of finite rank
of $\TT$. 

In particular, if $\und{u}_1,\ldots,\und{u}_m$ are elements of $\TT$, then there
exists an $S$-group $\Gamma$ of rank $n\leq m$ such that:
\[\und{u}_1^S\cdots\und{u}_m^S:=\{\und{u}_1^{\beta_1}\cdots\und{u}_m^{\beta_m},
\beta_1,\ldots,\beta_m\in S\}\subset\Gamma.\]
\label{lemme:sgroups}\end{Lemme}
{\bf Proof.} We prove the first part. Let $\{\und{u}_1,\ldots,\und{u}_m\}$ be a set of generators of $\Gamma$, let
$\Gamma_k$ be the $S$-submodule of $\Gamma$ generated by $\{\und{u}_1,\ldots,\und{u}_k\}$.
The $S$-module $\Gamma_1$ is clearly an $S$-group of finite rank (equal to $0$ if $\und{u}_1$ is
a torsion point of $\TT$, equal to $1$ if $\und{u}_1$ is a point of infinite order).

We now suppose that $\Gamma_{k-1}$ is contained in an $S$-group $\Lambda_{k-1}\subset\TT$ 
of rank $r\geq 0$ and we proceed to prove by induction that $\Gamma_k$ is also contained in some $S$-group
$\Lambda_{k}\subset\TT$ of rank $l\geq r$.

It is enough to construct an $S$-group of finite rank $\Lambda_{k}\supset\Lambda_{k-1}$ such that $\und{u}_k\in\Lambda_k$.
By hypothesis: \[\Lambda_{k-1}=(\Lambda_{k-1})_{{\tiny\mbox{ tors}}}\und{v}_1^S\cdots\und{v}_r^S,\]
where $(\Lambda_{k-1})_{{\tiny\mbox{ tors}}}$ is a finite torsion subgroup of $\TT$
and $\und{v}_1,\ldots,\und{v}_r$ are multiplicatively independent elements of $\TT$.

If $\und{v}_1,\ldots,\und{v}_r,\und{u}_k$ are multiplicatively independent, we put
\[\Lambda_k=\Lambda_{k-1}\und{u}_k^S.\] 
Let us suppose now that $\und{v}_1,\ldots,\und{v}_r,\und{u}_k$ are multiplicatively dependent. We have
a relation:\[\und{u}_k^\beta=\und{v}_1^{\beta_1}\cdots\und{v}_r^{\beta_r},\] where $\beta\in S-\{0\}$
and $(\beta_1,\ldots,\beta_r)\in S^r-\{(0,\ldots,0)\}$. 

For $i=1,\ldots,r$, let $\und{v}_i^\sharp$ be an element of $\TT$ such that:
\[\und{v}_i^\sharp{}^\beta=\und{v}_i.\]
Let $\Lambda_{k-1}^\sharp$ be the 
$S$-submodule of $\TT$ generated by $(\Lambda_{k-1})_{{\tiny\mbox{ tors}}}$, 
the kernel $\Ker(\beta)$ of the
isogeny given by $\und{u}\mapsto\und{u}^\beta$, and
the elements $\und{v}_i^\sharp$ for $i=1,\ldots,r$.

Since the $\und{v}_i^\sharp$ are multiplicatively independent, $\Lambda_{k-1}^\sharp$
is an $S$-group of rank $r$: it contains $\Lambda_{k-1}$ and $\und{u}_k$. We put 
$\Lambda_{k}=\Lambda_{k-1}^\sharp$. Since the second part of the lemma 
is clearly a consequence of the first part, the proof of the lemma is complete.

\subsubsection{Basic facts about $S$-groups and multiplicative independence.}

\begin{Lemme}
Let $\und{a}_1,\ldots,\und{a}_n$ be points of $\TT$,
let $(\xi_1,\xi_1'),\ldots,(\xi_n,\xi_n')$ be elements of $\CC^2$ such that 
$\Phi_0(\xi_i,\xi_i')=\und{a}_i$ for $i=1,\ldots,n$.
Then the following three conditions are equivalent.
\begin{enumerate}
\item The points $\und{a}_1,\ldots,\und{a}_n$ are multiplicatively dependent.
\item There exist elements $\tau_1,\ldots,\tau_n\in K$ not all zero and $\tau\in K$ such that, in $\CC$:
\begin{equation}\sum_{i=1}^n\tau_i\xi_i=\tau,\quad\mbox{ and }\sum_{i=1}^n\tau_i'\xi_i'=\tau'.\label{eq:relation_over_K}
\end{equation}
\item If $G$ is a finite torsion subgroup of $\TT$, then
$G\und{a}_1^S\cdots\und{a}_n^S$ is an $S$-group of rank $n$.
\item there exists an $n$-tuple
$\und{\tau}=(\tau_1,\ldots,\tau_n)\in K^n-\{\und{0}\}$ and 
two elements $\nu_1,\nu_2\in K$ linearly independent over $\QQ$,
such that:
\begin{equation}
\tr(\nu_1(\und{\tau}\cdot{}^{{\rm t}}\und{\xi})),
\quad\tr(\nu_2(\und{\tau}\cdot{}^{{\rm t}}\und{\xi}))\in\QQ,\label{eq:trtr}
\end{equation} where $\und{\tau}\cdot{}^{{\rm t}}\und{\xi}$ stands for
$\sum_{l=1}^n\tau_i\xi_i$.
\end{enumerate}
\label{lemme:independance_elementaire}
\end{Lemme}
{\bf Proof.} The equivalence between the first three conditions of the lemma is easily checked; for example,
the first two conditions are equivalent because 
for all $\tau\in K$ there exists $p\in\ZZ\setminus\{0\}$ such that $p\tau\in S$. 

We prove that
the fourth condition implies the second condition.
Let $\nu_1,\nu_2$ be two elements of $K$ which are
$\QQ$-linearly independent and $\und{\tau}\in K^n\setminus\{0\}$, satisfying (\ref{eq:trtr}).
 Let us suppose that:
\[\tr(\nu_1(\und{\tau}\cdot{}^{{\rm t}}\und{\xi}))=r_1,\quad
\tr(\nu_2(\und{\tau}\cdot{}^{{\rm t}}\und{\xi}))=r_2,\]
with $r_1,r_2\in\QQ$.
It is enough to prove that
$(\und{\mu}\cdot{}^{{\rm t}}\und{\xi},\und{\mu}'
\cdot{}^{{\rm t}}\und{\xi}')\in\Sigma(K)$.

Let us write $\xi_i=x_i+{\rm i}y_i$ with $x_i,y_i\in\RR$ for $i=1,\ldots,n$.
Let $\alpha=\und{\tau}\cdot{}^{{\rm t}}\und{x}$, 
$\alpha'=\und{\tau}'\cdot{}^{{\rm t}}\und{x}'$,
$\beta=\und{\tau}\cdot{}^{{\rm t}}\und{y}$ and $\beta'=\und{\tau}'\cdot{}^{{\rm t}}\und{y}'$. 
As the matrix $\sqm{\nu_1}{\nu_1'}{\nu_2}{\nu_2'}$ is invertible,
from
\[\sqm{\nu_1}{\nu_1'}{\nu_2}{\nu_2'}\cdot\binomial{\alpha+{\rm i}\beta}{\alpha'+{\rm i}\beta'}=\binomial{r_1}{r_2},\]
and separating real and imaginary part, we get $\beta=\beta'=0$ and:
\begin{eqnarray*}
\left({\alpha \atop \alpha'}\right) & = & (\nu_1\nu_2'-\nu_2\nu_1')^{-1}
\sqm{\nu_2'}{-\nu_1'}{-\nu_2}{\nu_1}\cdot\left({r_1 \atop r_2}\right)\\
& = & \Sigma\left(\frac{\nu_2'r_1-\nu_1'r_2}{\nu_1\nu_2'-\nu_2\nu_1'}\right)\in\Sigma(K).
\end{eqnarray*}
This is the second condition. If on the other hand the relations
(\ref{eq:relation_over_K}) hold, then, we have
$(\und{\tau}\cdot{}^{\rm t}\und{\xi},\und{\tau}'\cdot{}^{\rm t}\und{\xi}')\in\Sigma(K)$, and
for any two elements $\nu_1,\nu_2$ of $K$ we have (\ref{eq:trtr}).
The lemma \ref{lemme:independance_elementaire} is proven.

\subsubsection{Algebraic interpretation of an analytic property.}

The lemma below is a generalisation of lemma 3.2 p. 212
of \cite{Masser:Hecke}.

\begin{Lemme}
Let us consider points $\und{a}_1,\ldots,\und{a}_n\in\TT(\bar{\QQ})\cap{\cal D}$.
Then, the point \[\und{A}:=(\und{a}_1,\ldots,\und{a}_n)\in\TT^n\] is in ${\cal U}({\cal B})$.

Moreover, the points $\und{a}_1,\ldots,\und{a}_n$ are multiplicatively independent 
if and only if $\und{A}$ satisfies the property A.
\label{lemme:proprieteA}\end{Lemme}
{\bf Proof.}
It is easy to see that ${\cal D}\subset{\cal U}({\cal B}_0(\eta))$: using the same arguments as in lemma 3.2 p. 212 of
\cite{Masser:Hecke}, we also see that $\und{A}\in{\cal U}({\cal B})$.

Let us suppose that $\und{A}$ does not
satisfy the property A. 
The theorem p. 276 of \cite{Masser:Vanishing} implies that
there exist two distinct elements \[(p_{1,1},q_{1,1},\ldots, p_{n,1},q_{n,1}),
(p_{1,2},q_{1,2},\ldots,p_{n,2},q_{n,2})\in\NN^{2n},\] 
and an arithmetic progression ${\goth K}\subset \NN$, such that if we write:
\[M_i(U_1,V_1,\ldots,U_n,V_n) =
U_1^{p_{1,i}}V_1^{q_{1,i}}\cdots U_n^{p_{n,i}}V_n^{q_{n,i}} \mbox{ for $i=1,2$ },\]
then we have:
\begin{small}\begin{equation}
M_1({\cal B}^\kappa.(u_1,v_1,\ldots,u_n,v_n))=M_2({\cal B}^\kappa.(u_1,v_1,\ldots,u_n,v_n))
\mbox{ for all $\kappa\in{\goth K}$}.\label{eq:premier}\end{equation}\end{small}

Let $(\xi_1,\xi_1'),\ldots,(\xi_n,\xi_n')\in\CC^2$ be such that
$\Phi_0(\xi_i,\xi_i')=\und{a}_i=(u_i,v_i)$, let
$\und{\mu}_1,\und{\mu}_2\in M^*{}^{n}$ such that
\[\Sigma^{\oplus n}(\und{\mu}_i)=(p_{1,i},q_{1,i},\ldots,p_{n,i},q_{n,i})\cdot{\goth B}_0^{\oplus n}
\mbox{ for $i=1,2$}.\] 

We observe that $\und{\tau}=\und{\mu}_1-\und{\mu}_2\not=\und{0}$.
We have~:
\begin{eqnarray*}
\lefteqn{M_i({\cal B}^\kappa.(u_1,v_1,\ldots,u_n,v_n))=}\\
& = &\exp\{2\pi{\rm i}(p_{1,i},q_{1,i},\ldots,p_{n,i},q_{n,i})\cdot {\cal B}^\kappa
\cdot{\goth B}_0^{\oplus n}\cdot{}^{{\rm t}}(\xi_1,\xi_1',\ldots,\xi_n,\xi_n')\}\\
& = &\exp\{2\pi{\rm i}(p_{1,i},q_{1,i},\ldots,p_{n,i},q_{n,i})\cdot
{\goth B}_0^{\oplus n}\cdot{}^{{\rm t}}(\eta^{\kappa}\xi_1,\eta^{-\kappa}\xi_1',\ldots,
\eta^{\kappa}\xi_n,\eta^{-\kappa}\xi_n')\}\\
& = & \exp\{2\pi{\rm i}\tr(\eta^{\kappa}\und{\mu}_i\cdot{}^{{\rm t}}\und{\xi})\}
\mbox{ for $i=1,2$.}\end{eqnarray*}
The identities (\ref{eq:premier}) for $\kappa\in{\goth K}$ imply:
\[\tr(\eta^{\kappa_1}\und{\tau}\cdot{}^{{\rm t}}\und{\xi}),
\quad\tr(\eta^{\kappa_2}\und{\tau}\cdot{}^{{\rm t}}\und{\xi})\in\QQ\]
for two distinct rational integers $\kappa_1$ and $\kappa_2$. But 
$\eta^{\kappa_1}$ and $\eta^{\kappa_2}$ are $\QQ$-linearly 
independent and lemma \ref{lemme:independance_elementaire} implies 
that $\und{a}_1,\ldots,\und{a}_n$ are multiplicatively dependent.

On the other side, if $\und{a}_1,\ldots,\und{a}_n\in{\cal D}$ are multiplicatively dependent,
then it is easy to see that for some $l\in\NN^\times$, the point 
\[\und{A}'=\und{A}^l=(\und{a}_1^l,\ldots,\und{a}_n^l)\] 
lies in a connected algebraic
subgroup $H$ of $\TT^n\cong\GG_m(\CC)^{2n}$ of even codimension $\geq 2$, stable under the 
diagonal action of $S$ (over $\TT^n$). In particular, for all $k\in\NN$,
\[\und{A}'{}^{\und{\eta}^k}\in H.\] Let $H_1$ be a hypersurface
containing $H$; for all $k$, \[{\cal B}^k.\und{A}'=\und{A}'{}^{\und{\eta}^k}\in H_1.\]
The reader can check that $H_1$ can be choosen so that an equation defining it also provides
a non zero locally analytic function $F$ such that for all $k$ big enough,
\[F({\cal B}^k.\und{A}')=0.\] This condition implies that the point $\und{A}$ does not satisfy 
the property $A$.

\medskip

\noindent {\bf Remark.} In this text we only need to know that if $\und{a}_1,\ldots,\und{a}_n$
are multiplicatively independent, then the point $\und{A}$ satisfies the property $A$.

\section{End of proof of theorem 2.}

Let $\und{u}_1,\ldots,\und{u}_m$ be elements of infinite order of $\TT(\bar{\QQ})\cap {\cal D}$,
let us consider the $S$-submodule $\Lambda$ of $\TT$ generated by 
$\und{u}_1,\ldots,\und{u}_m$ and let us choose an $S$-group of finite rank $\Gamma$ of $\TT(\bar{\QQ})$ 
such that $\Gamma\supset\Lambda$, whose existence is guaranteed by lemma \ref{lemme:sgroups}:
let $n\geq 1$ be its rank.
Let $W$ be any infinite group of units $\tilde{\eta}$ of $S_+$
acting as the identity map on
$\Gamma_{\tiny {\mbox{{\rm tors}}}}$ (it exists thanks to lemma \ref{lemma:torsion_3}). 

\subsection{Looking for a positive $S$-basis.}

We will need the following proposition, which is a variant for a real multiplication structure, of 
lemma 3 p. 37 of \cite{Loxton:Fredholm}.
\begin{Proposition}
It is always possible to find multiplicatively independent elements
$\und{a}_1,\ldots,\und{a}_n\in\TT(\bar{\QQ})\cap{\cal D}^+$, elements 
$\und{\nu}_j=(\nu_{j,1},\ldots,\nu_{j,n})$ such that $\Sigma(\nu_{i,j})\in {\cal X}$,
elements $\eta_1,\ldots,\eta_m\in W$, 
and elements $\alpha_j\in K$ such that:
\begin{equation}
\und{u}_i^{\eta_i}=\Phi_0(\alpha_i,\alpha_i')\und{a}_1^{\nu_{i,1}}\cdots\und{a}_n^{\nu_{i,n}},\quad
\label{eq:ritt}\mbox{ for }i=1,\ldots,m.\end{equation}\label{lemme:ritt}\end{Proposition}
The proof of lemma 3 p. 37 of \cite{Loxton:Fredholm} is simple. Conversely, our proof
of proposition \ref{lemme:ritt} is quite involved (Baker's theorem
on linear forms of two logarithms of algebraic numbers occurs), and it is more convenient to begin with some 
technical settings, and to divide the proof of the proposition in several lemmata.

\subsubsection{Topology of matrices.}

We fix an isomorphism between the euclidean vector spaces $\RR^{n^2}$ and $M_{n\times n}(\RR)$,
the space of square matrices of order $n$ and real entries.

\begin{Lemme}
Let $\und{\alpha},\und{y}\in\RR^n$ be non zero column matrices, let $T_0$ be a matrix of
$M_{n\times n}(\RR)$ such that \[T_0\cdot\und{\alpha}=\und{y}.\] There exists a non empty subset
\[{\cal A}_{\und{\alpha},\und{y}}\subset M_{n\times n}(\RR)\setminus Z,\]
locally isomorphic to $\RR^{n^2-n}$, such that $T_0\in\overline{{\cal A}_{\und{\alpha},\und{y}}}$ (euclidean adherence),
and such that for all $T\in{\cal A}_{\und{\alpha},\und{y}}$,
$T\cdot\und{\alpha}=\und{y}$.
\label{lemme:foglio1}\end{Lemme}

\noindent {\bf Proof.} Since $\und{\alpha},\und{y}$ do not vanish, the set
\[{\cal T}_{\und{\alpha},\und{y}}=\{T\in M_{n\times n}(\RR)\mbox{ such that }T\cdot\und{\alpha}=\und{y}\}\]
is a proper translated of a vector subspace of $M_{n\times n}(\RR)$ of dimension $n^2-n$.

Clearly, ${\cal T}_{\und{\alpha},\und{y}}\not\subset Z$,
where $Z$ is the locus of vanishing of the determinant in $M_{n\times n}(\RR)$.
Thus, 
\[{\cal A}_{\und{\alpha},\und{y}}={\cal T}_{\und{\alpha},\und{y}}\setminus Z\] is non-empty and locally 
isomorphic to $\RR^{n^2-n}$. 

If $\det(T_0)\not=0$, then clearly, $T_0\in{\cal A}_{\und{\alpha},\und{y}}$.
It $\det(T_0)=0$ it is easy to construct a sequence $(T_i)_{i\geq 1}$ of matrices $T_i\in
{\cal A}_{\und{\alpha},\und{y}}$, such that $\lim_{i\rightarrow\infty}T_i=T_0$.

\medskip

In an euclidean vector space $V$, we denote $B(\und{x},r)$ the open euclidean ball of center 
$\und{x}$ and radius $r>0$.

\begin{Lemme}
Let ${\cal O}$ be a non-empty open subset of $\RR^n$, $\und{\alpha},\und{y}\in\RR^n$ be non zero
wit $\und{\alpha}\in\overline{{\cal O}}$, let $T_0\in M_{n\times n}(\RR)$
be such that \[T_0\cdot\und{\alpha}=\und{y}.\]

For $\epsilon>0$, let ${\cal A}_{\und{\alpha},\und{y},\epsilon}$ be the subset 
of $\RR^n$ whose elements are the matrices \[T\in (M_{n\times
n}(\RR)\setminus Z)\cap B(T_0,\epsilon)\] such that there exists $\und{\alpha}^*\in B(\und{\alpha},\epsilon)\cap{\cal O}$ with
$T\cdot\und{\alpha}^*=\und{y}$.

For all $\epsilon$ small enough, ${\cal A}_{\und{\alpha},\und{y},\epsilon}$ is non-empty, open, and satisfies the
property that $T_0\in\overline{{\cal A}_{\und{\alpha},\und{y},\epsilon}}$.
\label{lemme:foglio2}\end{Lemme}
\noindent {\bf Proof.} For $\epsilon>0$ small enough, $\und{0}\not\in B(\und{\alpha},\epsilon)$.
Let $\und{\alpha}^*\in B(\und{\alpha},\epsilon)\cap{\cal O}$; then $\und{\alpha}^*\not=0$ and we can apply lemma \ref{lemme:foglio1}.
Let 
${\cal A}_{\und{\alpha}^*,\und{y}}$ be the set given by lemma \ref{lemme:foglio1}. We see from the definition
(the elements of ${\cal A}_{\und{\alpha},\und{y}}$ are invertible matrices), that
if $\und{\alpha_1}\not=\und{\alpha_2}$ are two elements of $B(\und{\alpha},\epsilon)$ or $B(\und{\alpha},\epsilon)\cap{\cal O}$, then:
\[{\cal A}_{\und{\alpha}_1,\und{y}}\cap{\cal A}_{\und{\alpha}_2,\und{y}}=\emptyset.\]
As $B(\und{\alpha},\epsilon)\cap{\cal O}$ is locally isomorphic to $\RR^n$, the set 
\[\tilde{{\cal A}}=\bigcup_{\und{\alpha}^*\in B(\und{\alpha},\epsilon)\cap{\cal O}}{\cal A}_{\und{\alpha}^*,\und{y}}\quad \mbox{(disjoint
union)},\] is locally isomorphic to $\RR^{n^2}$. 
Moreover, \[{{\cal A}_{\und{\alpha},\und{y},\epsilon}}=\tilde{{\cal A}}\cap
B(T_0,\epsilon),\] thus non-empty and open.
From lemma \ref{lemme:foglio1} we also see that $T_0\in\overline{\und{\cal A}_{{\alpha},\und{y},\epsilon}}$.

\subsubsection{Two open sets and their intersection.}

Let $T,T'\in M_{n\times n}(\RR)$:
\begin{equation}T=\left(\begin{array}{ccc}\tau_{1,1} & \cdots & \tau_{1,n}\\
\vdots  & & \vdots \\
\tau_{n,1}  & \cdots & \tau_{n,n}
\end{array}\right),\quad T'=\left(\begin{array}{ccc}\tau_{1,1}' & \cdots & \tau_{1,n}'\\
\vdots &  & \vdots \\
\tau_{n,1}' &  \cdots & \tau_{n,n}'
\end{array}\right).\label{eq:matrices}\end{equation}
Then we define:
\[L_{T,T'}=\left(\begin{array}{ccc}b_{1,1} & \cdots & b_{1,n}\\
\vdots &  & \vdots \\
b_{n,1} & \cdots & b_{n,n}
\end{array}\right)\in M_{2n\times 2n}(\RR),\] where the block $b_{i,j}$ is the matrix 
$\sqm{\tau_{i,j}}{0}{0}{\tau_{i,j}'}\in M_{2\times 2}(\RR)$.

Let
\[\lambda_{T,T'}:\RR^{2n}\rightarrow\RR^{2n}\] be the linear map 
associated to the matrix $L_{T,T'}$, acting on column vectors from the left. The map $\lambda_{T,T'}$
is an automorphism if and only
if both of
$T,T'$ are non-singular, because $\det(L_{T,T'})=\det(T)\det(T')$.

Let $\delta$ be a positive number, let us choose: 
\begin{eqnarray}
\und{y}&=&(y_1,\ldots,y_n),\in\RR^n\setminus\{\und{0}\}\nonumber\\
\und{y}'&=&(y_1',\ldots,y_n'),\in\RR^n\setminus\{\und{0}\}\label{eq:yY}\\
\und{Y}& = & (y_1,y_1',\ldots,y_n,y_n')\in(\RR^n\setminus\{\und{0}\})^2,\nonumber
\end{eqnarray}
and let us define the matrices
\[T_0=\left(\begin{array}{lll}y_1 & \cdots & y_1\\
\vdots & & \vdots\\ y_n & \cdots & y_n\end{array}\right)\mbox{ and }
T_0'=\left(\begin{array}{lll}y_1' & \cdots & y_1'\\
\vdots & & \vdots\\ y_n' & \cdots & y_n'\end{array}\right).\] 
We introduce two subsets ${\goth X}(\delta),{\goth Y}$ of 
the real euclidean space $M_{n\times n}(\RR)\times M_{n\times n}(\RR)\cong\RR^{2n^2}$ of couples of square matrices
$(T,T')$ of order $n$.

\vspace{10pt}

\noindent ${\goth X}(\delta)$ is the non-empty open set of couples $(T,T')$ with $T=(\tau_{i,j}),T'=(\tau_{i,j}')$
such that the absolute values of the
coefficients of the matrices:
\[T-T_0,\quad T'-T_0'\] are $<\delta$.

\vspace{10pt}

\noindent ${\goth Y}$ is the set of all the couples of regular matrices $(T,T')$ such that 
\[\und{Y}\in\lambda_{T,T'}(\Im({\cal W}^+)^n),\] where $\Im({\cal W}^+)=
\{(\Im(z),\Im(z'))\mbox{ with }(z,z')\in{\cal W}^+\}\subset\RR^2$. 

\begin{Lemme} 
The set ${\goth Y}$ is a non-empty open set such that the couple of matrices
$(T_0,T'_0)$ belongs to its 
euclidean adherence $\overline{{\goth Y}}$. 

For all $\delta>0$,
the intersection 
\begin{equation}{\goth Z}(\delta):={\goth X}(\delta)\cap {\goth Y}
\label{eq:zeta}\end{equation}
contains a non-empty open set.
\label{lemme:topology}\end{Lemme}

\noindent {\bf Proof.}
We prove that for all $\epsilon>0$ small enough,
the set
\[B((T_0,T_0'),\epsilon)\cap{\goth Y}\]
contains a non-empty open subset. 

Indeed, let us choose
\[\und{\alpha}=\und{\alpha}'= \left(\frac{1}{n},\ldots,\frac{1}{n}\right),\]
and let us denote:
\begin{eqnarray*}
{\cal O} & = & \{\und{x}=(x_1,\ldots,x_n)\mbox{ such that }x_i>1/n\},\\
{\cal O}' & = & \{\und{x}'=(x_1',\ldots,x_n')\mbox{ such that }-1/n<x_i'<1/n\}.
\end{eqnarray*}
We apply lemma \ref{lemme:foglio2} twice: once for
the data $T_0,{\cal O},\und{\alpha},\und{y}$, once for the data $T_0',{\cal O}',\und{\alpha}',\und{y}'$ 
(all the hypotheses are satisfied):
we obtain, for $\epsilon>0$ small enough, that
\[{\cal A}':={\cal A}_{\und{\alpha},\und{y},\epsilon}\times{\cal A}_{\und{\alpha}',\und{y}',\epsilon}\]
is non-empty and open in $\RR^{2n^2}$, and the reader can easily check that it is contained in ${\goth Y}\cap
B((T_0,T_0'),\epsilon)$.

The set ${\goth Y}$ is open: let 
$(T,T')\in{\goth Y}$. If ${\cal E},{\cal E}'$ are matrices of size $n\times n$ whose
coefficients are real numbers of small enough absolute values, then $T+{\cal E}$ and $T'+{\cal E}'$
are regular, because the locus of couples of matrices $(S,S')$ such that $\det(S)\det(S')=0$
determines a Zariski proper closed subset of $\RR^{2n^2}$ not containing $(T,T')$. Let $C(\epsilon)$ be a little
cube of size
$\epsilon>0$ centered at
$\und{0}\in\RR^{2n^2}$.

If $\epsilon>0$ is small enough, the function $H:\RR^{2n^2}\rightarrow\RR^{2n}$
defined by:
\[H:({\cal E},{\cal E}')\mapsto(\lambda_{T+{\cal E},T'+{\cal E}'})^{-1}(\und{Y})\] is well defined and of class
${\cal C}^\infty$ on $C(\epsilon)$. Since $H(\und{0},\und{0})\in\Im({\cal W}^+)^n$ by hypothesis, and since 
the latter set is open and non-empty, for $\epsilon$ small enough $H(C(\epsilon))\subset\Im({\cal W}^+)^n$
and $C(\epsilon)\subset{\goth Y}$. 

Since for all $\delta>0$, $(T_0,T_0')\in{\goth X}(\delta)$, we see that for all $\delta>0$, 
the intersection ${\goth Z}(\delta)$
 contains 
a non-empty open set. The lemma is proved. 

\subsubsection{Application of Baker's theorem.}
\begin{Lemme} Let $\und{b}=(a,b)\in\TT(\bar{\QQ})$; then the following conditions are equivalent.
\begin{enumerate}
\item $|a||b|^\theta=1$.
\item $|a||b|^{\theta'}=1$.
\item $|a|=|b|=1$.
\end{enumerate}
\label{lemme:baker}\end{Lemme}
{\bf Proof.} Since $|a|,|b|$ are algebraic numbers, If $|a|\not=1$ or $|b|\not=1$, then $\log|a|+\theta\log|b|=0$
if and only if there is a rational number $r$ such that the matrix $\displaystyle{\sqm{1}{\theta}{1}{r}}$
is singular, by using the ineffective Baker's theorem on linear forms of logarithms
of algebraic numbers (see \cite{Baker:Logarithms} for the foundations of the theory), because $|a|,|b|$ are real algebraic
numbers. But $\theta$ is irrational, and the matrix above cannot be singular. In this way we see that 
the first and the third conditions of the lemma are equivalent.

The same technique can be applied to the linear form $\log|a|+\theta'\log|b|$, to prove that also the 
second and the third conditions of the lemma are equivalent.

\subsubsection{Proof of proposition \ref{lemme:ritt}.}

It is easy to check, by using lemma \ref{lemma:set_S}, that:
\[{\cal X}\subset\Im({\cal W}^+)\subset\Im({\cal W}).\] 
Since $\und{u}_i\in{\cal D}$ for $i=1,\ldots,m$,
for all $i=1,\ldots,m$ there exists $\eta_i\in W$ such that
$\und{u}_i^{\eta_i}=\Phi_0(w_i,w_i')$, for a couple of complex numbers $(w_i,w_i')$ with 
\[(t_i,t_i'):=(\Im(w_i),\Im(w_i'))\in{\cal X}.\] 

The $S$-module generated by the elements 
$\und{u}_i^{\eta_i}$ is also contained in $\Gamma$ because of the choice of $W$, which
acts as the identity on the torsion subgroup. Let us choose an $S$-basis $(\und{b}_1,\ldots,\und{b}_n)$
of $\Gamma/\Gamma_{{\tiny \mbox{tors}}}$, i. e. a maximal collection of multiplicatively independent 
elements of $\Gamma$.

There exist elements $\alpha_1,\ldots,\alpha_m\in K$ and elements $\und{\mu}_1,\ldots,\und{\mu}_m
\in S^n-\{\und{0}\}$ with $\und{\mu}_i=(\mu_{i,1},\ldots,\mu_{i,n})$ such that:
\begin{equation}
\und{u}_i^{\eta_i}=\Phi_0(\alpha_i,\alpha_i')\und{b}_1^{\mu_{i,1}}\cdots
\und{b}_n^{\mu_{i,n}}.\label{eq:b}\end{equation}
Let us fix elements $(z_1,z_1'),\ldots,(z_n,z_n')\in\CC^2$ such that 
\[\Phi_0(z_i,z_i')=\und{b}_i=(a_i,b_i)\in\TT(\bar{\QQ}).\] 
Let us write $y_i=\Im(z_i)$, $y_i'=\Im(z_i')$; the identities (\ref{eq:b}) imply:
\begin{equation}
t_i=\sum_{j=1}^n\mu_{i,j}y_j,\quad t_i'=\sum_{j=1}^n\mu_{i,j}'y_j'.
\label{eq:riscrittura}\end{equation}
Let us consider $\und{y},\und{y}'$ and $\und{Y}$ as in (\ref{eq:yY}), let \[\und{t}=(t_1,t_1'\ldots,t_m,t_m')\in\RR^{2m}\]
be a column vector, let 
\[M=\left(\begin{array}{ccc}\mu_{1,1} & \cdots & \mu_{1,n}\\
\vdots & & \vdots\\
\mu_{m,1} & \cdots & \mu_{m,n}\end{array}\right)\in M_{m\times n}(S).\]
We may rewrite (\ref{eq:riscrittura}) as follows:
\begin{eqnarray}\und{t} & = & L_{M,M'}\cdot\und{Y}\nonumber\\
& = & (L_{M,M'}\cdot L_{T,T'})\cdot(L_{T,T'}^{-1}\cdot\und{Y})\label{eq:matrice_T}
\end{eqnarray}
for every matrix $T\in\GL_n(K)$ (in this case, $T'$ is its Galois conjuguate). We now proceed to choose $T$ 
in a fruitful way:
let the coefficients of $T$ and $T'$ be as in (\ref{eq:matrices}).

We see that both $\und{y},\und{y}'$ are non zero.
Indeed, applying lemma \ref{lemme:baker}, $y_i=0$ (for some $i$) if and
only if $|a_i|=|b_i|=1$ if and only if $y_i'=0$. Thus $\und{y}=\und{0}$ if and only if $\und{y}'=\und{0}$.
But if $\und{y},\und{y}'$ both vanish, then 
the $S$-group generated by $\und{b}_1,\ldots,\und{b}_n$ would be entirely
contained in the boundary of ${\cal D}$, and in this case, it could not contain elements $\und{u}_i\in{\cal D}$
(\footnote{Note however, that some elements of the $S$-basis $(\und{b}_1,\ldots,\und{b}_n)$ might belong to the
boundary of ${\cal D}$.}):
lemma \ref{lemme:topology} can be applied, and the set ${\goth Z}(\delta)$ defined in (\ref{eq:zeta})
is non-empty and open, for all $\delta>0$.

By lemma \ref{lemme:topology}, 
${\goth Z}(\delta)$ contains a non-empty open set.
Hence, for all $\delta>0$ 
there exists $T\in\GL_n(K)$ such that $(T,T')\in{\goth Z}(\delta)$, because
$\Sigma(K)$ is dense in $\RR^2$ for the euclidean topology. We fix such a couple of matrices.

First of all, since $(T,T')\in{\goth Y}$, we have:
\begin{equation}\und{\tilde{Y}}:=L_{T,T'}^{-1}\cdot\und{Y}\in\Im({\cal W}^+)^n.\label{eq:gothy}\end{equation}

Moreover, if $\delta>0$ is small enough depending only
on the elements $\und{\mu}_i$ and on the basis $\und{b}_1,\ldots,\und{b}_n$, then the coefficients $\varsigma_{i,j}$
of the matrix $M\cdot T=(\varsigma_{i,j})_{i,j}$ satisfy, by using $L_{M,M'}\cdot L_{T,T'}=L_{M\cdot T,M'\cdot T'}$:
\begin{equation}\Sigma(\varsigma_{i,j})\in{\cal X}\quad\mbox{ for all }i,j.\label{eq:varsigma_i}\end{equation}
In effect, since we have the identities (\ref{eq:riscrittura}) in $\RR$ and (in $K$):
\[\varsigma_{i,j}=\sum_{k=1}^n\mu_{i,k}\tau_{k,j},\] if $(T,T')$ and $(T_0,T_0')$
are close enough (if $\delta$ is small enough), then $\Sigma(\varsigma_{i,j})$ is close enough to
$(t_i,t_i')$ ($i=1,\ldots,m,j=1,\ldots,m$) to lie in ${\cal X}$, since $(t_i,t_i')\in{\cal X}$ 
($i=1,\ldots,m$) by hypothesis.

Hence, (\ref{eq:gothy}) and (\ref{eq:varsigma_i}) hold at once for our choice of $T\in\GL_n(K)$.

We now consider $\lambda_{T,T'}$ as an automorphism $\CC^{2n}\rightarrow\CC^{2n}$. 
Let us denote (column matrices): \[(\tilde{z}_1,\tilde{z}'_1,\ldots,\tilde{z}_n,\tilde{z}'_n)=\lambda_{T,T'}^{-1}
(z_1,z_1',\ldots, z_n,z_n').\] Since $(z_1,z_1',\ldots, z_n,z_n')\in\lambda_{T,T'}(({\cal W}^+)^n)$,
we have that $(\tilde{z}_i,\tilde{z}_i')\in
{\cal W}^+$ for all $i=1,\ldots,n$, because of (\ref{eq:gothy}). 

From (\ref{eq:varsigma_i}) we see that:
\[\und{\varsigma}_i:=(\varsigma_{i,1},\varsigma_{i,1}',\ldots,\varsigma_{i,n},\varsigma_{i,n}')
\in({\cal X}\cap\Sigma(K))^n\quad (i=1,\ldots,m).\]

A Kronecker-type extension provides the end of the proof.
Let $d\in\ZZ_{>0}$ be such that $d\varsigma_{i,j}\in S$ for all $i,j$ (in this way, $\Sigma(d\varsigma_{i,j})\in{\cal X}$
because $\Sigma(\varsigma_{i,j})\in{\cal X}$),
let  $\xi_i=\tilde{z}_i/d$, $\xi_i'=\tilde{z}_i'/d$ (we see that $(\xi_i,\xi_i')\in{\cal W}^+$).
The proof of the proposition \ref{lemme:ritt} is complete by setting~:
\[\und{a}_j=\Phi_0(\xi_j,\xi_j')\in{\cal D}^+,\nu_{i,j}=d\varsigma_{i,j}.\]
Indeed, for $i=1,\ldots,m$:
\begin{eqnarray*}
\und{u}_i^{\eta_i} & = & \Phi_0\left(\alpha_i+\sum_{j=1}^n\varsigma_{i,j}\tilde{z}_j,
\alpha_i'+\sum_{j=1}^n\varsigma_{i,j}'\tilde{z}_j'\right)\\
& = & \Phi_0(\alpha_i,\alpha_i')\Phi_0\left(\sum_{j=1}^n\nu_{i,j}\frac{\tilde{z}_j}{d},
\sum_{j=1}^n\nu_{i,j}'\frac{\tilde{z}_j'}{d}\right)\\
& = & \Phi_0(\alpha_i,\alpha_i')\und{a}_1^{\nu_{i,1}}\cdots\und{a}_n^{\nu_{i,n}}.
\end{eqnarray*}


\subsection{Auxiliary functions.\label{section:auxiliary}}
We are ready to continue and complete the proof of the theorem 2.
Let $\und{u}_1,\ldots,\und{u}_m$ be elements of $\TT(\bar{\QQ})\cap{\cal D}$,
let us suppose by contradiction 
that the complex numbers $f(\und{u}_1),\ldots,f(\und{u}_m)$ are algebraically dependent over $\QQ$. 

In view of an application of the proposition \ref{propo:loxton} of the section \ref{section:tools}, we must exhibit a 
certain choice of locally analytic functions $\Psi_i$ for $i=1,\ldots,m$.

We proceed like this. We apply the proposition \ref{lemme:ritt}:
let $W$ be as in the proposition. We then have a positive integer $n$, elements
$\und{a}_1,\ldots,\und{a}_n\in\TT(\bar{\QQ})\cap
{\cal D}^+$ which are multiplicatively independent, elements $\und{\nu}_j=(\nu_{j,1},\ldots,\nu_{j,n})$
with $\Sigma(\nu_{j,k})\in{\cal X}$, 
elements
$\eta_1,\ldots,\eta_m\in W$ and elements $\alpha_j\in K$ so that the equalities 
(\ref{eq:ritt}) hold, with $\und{\zeta}_i=\Phi_0(\alpha_i,\alpha'_i)$.

We may arrange the indexes $i=1,\ldots,m$ so that there exists $0\leq m_0\leq m$ with:
\begin{eqnarray}
\nu_{i,1}&\in & S_+\quad\mbox{ for }i=1,\ldots,m_0,\label{eq:conditions_nu}\\
\nu_{i,1}&\in & S_\pm\mbox{ for }i=m_0+1,\ldots,m.\nonumber
\end{eqnarray}
According with the value of $m_0$, we define:
\begin{eqnarray}
\Psi_i(\und{V})= \Psi_i(\und{v}_1,\ldots,\und{v}_n) & = & f^+(\und{\zeta}_i\und{v}_1^{\nu_{i,1}}
\cdots\und{v}_n^{\nu_{i,n}})\mbox{ for $i=1,\ldots,m_0$,}\label{eq:psi}\\
& = & f(\und{\zeta}_i\und{v}_1^{\nu_{i,1}}\cdots\und{v}_n^{\nu_{i,n}})\mbox{ for $i=m_0+1,\ldots,m$.}
\nonumber\end{eqnarray}
Let us now check that the proposition \ref{propo:loxton} can be applied.

First of all, as $f(\und{u}),f^+(\und{u})$ are locally analytic and $\Sigma^{\oplus n}(\und{\nu}_i)\in{\cal X}^n$,
the functions $\Psi_i(\und{V})$ are locally analytic, thanks to lemma \ref{lemme:action_S}.

It is easy to see that the Taylor coefficients of the series defining the functions $\Psi_i(\und{V})$
are algebraic numbers in some cyclotomic number field of finite degree over $\QQ$.

The functions $\Psi_i(\und{V})$ converge simultaneously on a non-empty euclidean
open neighbourhood $\Omega$ of $\und{0}\in\CC^{2n}$ with the property that:
\[\Omega\subset{\cal U}({\cal B}).\]
After the identities (\ref{eq:functional_equation}) and
(\ref{eq:functional_equation2}), the functions $\Psi_i(\und{V})$ 
satisfy on $\Omega$ the collection of simultaneous functional equations:
\begin{equation}
\Psi_i(\und{V}^{\und{\eta}})=\Psi_i(\und{V})-R_{i}(\und{V}),\label{eq:funct_equations}
\end{equation}
where 
\begin{eqnarray*}
R_{i}(\und{V}) & = & R_{\eta}^+(\und{\zeta}_i\und{v}_1^{\nu_{i,1}}
\cdots\und{v}_n^{\nu_{i,n}})\mbox{ for $i=1,\ldots,m_0$,}\\
& = & R_{\eta}(\und{\zeta}_i\und{v}_1^{\nu_{i,1}}
\cdots\und{v}_n^{\nu_{i,n}})\mbox{ for $i=m_0+1,\ldots,m$.}
\end{eqnarray*}
Moreover, the matrix ${\cal B}$ involved
in the automorphism of (\ref{eq:funct_equations}) is good, as we said earlier, in section \ref{section:analytic_tools}.

Let us denote:\[\und{A}=(\und{a}_1,\ldots,\und{a}_n)\in\TT^n.\] By the first part of lemma 
\ref{lemme:proprieteA}, $\und{A}\in{\cal U}({\cal
B})$ so that for $k$ big enough, \[{\cal B}^k.\und{A}=\und{A}^{\und{\eta}^k}\in\Omega.\] 
Hence, for any $k$ big enough (as soon as the values make sense) we get:
\[\QQ(\Psi_1(\und{A}^{\und{\eta}^k}),\ldots,\Psi_m(\und{A}^{\und{\eta}^k}))=
\QQ(f(\und{u}_1),\ldots,f(\und{u}_m)),\] by an iterated application of the functional
equations (\ref{eq:funct_equations}), as well
as the identity $\Theta_M=f+f^+$.
As $\und{a}_1,\ldots,\und{a}_n$ are multiplicatively independent we see
that, after the second part of lemma \ref{lemme:proprieteA}, the point $\und{A}$ satisfies the
property A as well as $\und{A}^{\und{\eta}^k}$ for all $k>0$. Therefore, the proposition
\ref{propo:loxton} can be applied and says that the functions $\Psi_i(\und{V})$ are 
algebraically dependent over $\CC(\und{V})$.
\subsection{Algebraic independence of locally analytic functions.\label{section:algebraic}}
To end the proof of theorem 2 we still need to prove that, supposing the 
functions $\Psi_i(\und{V})$ defined in (\ref{eq:psi})
being algebraically dependent over $\CC(\und{V})$, the $m$-tuple ${\cal M}$
is not semi-free.

The corollary 9 p. 29 of \cite{Kubota:Algebraic} implies that the functions
$\Psi_i(\und{V})$ are $\CC$-{\em linearly dependent modulo} $\CC(\und{V})$ 
in the following sense. There exists $m$ complex numbers $c_1,\ldots,c_m$,
not all zero, and a rational function $Q(\und{V})\in\CC(\und{V})$ such that:
\begin{equation}
\sum_{i=1}^mc_i\Psi_i(\und{V})=Q(\und{V}).\label{eq:rel_dep_lin}
\end{equation} 
We still need to prove:

\begin{Proposition}
If the functions
$\Psi_i(\und{V})$ are $\CC$-{\em linearly dependent modulo} $\CC(\und{V})$, and satisfy (\ref{eq:rel_dep_lin}), then
the $m$-tuple ${\cal M}$ is not semi-free.
\label{proposition:linear_dep}\end{Proposition}

The plan of the proof of proposition \ref{proposition:linear_dep} is the following. We first prove that
the relations (\ref{eq:rel_dep_lin}) are equivalent
to $\CC$-linear dependence of rational functions $R_{\eta,N_i},R_{\eta,N_i}^+$.

Then we prove that the $\CC$-linear relation so obtained is defined over $\QQ$.

\medskip

\noindent {\bf Proof of proposition \ref{proposition:linear_dep}.\label{section:from_general_to_1}}
We clearly have, from (\ref{eq:conditions_nu}) and the definition (\ref{eq:psi}), the inclusions:
\begin{eqnarray}
\spt_K(\Psi_i) & \subset & \und{\nu}_{i}M^*\label{eq:support_Psi_K}\\
& \subset & (\nu_{1,i}M^*\cap{\cal I}_+)\cap (M^*)^{n-1}\nonumber\\
& \subset & (M^*)^n\cap{\goth H}_n^+,
\label{eq:support_Psi_i}\end{eqnarray}
for all
$i=1,\ldots,m$; moreover, there exists a strictly convex cone $\Pi\subset\RR_{\geq 0}\times\RR^{n-1}$
(of axis $\und{Y}=(1,0,\ldots,0)\RR$) such that for all $i=1,\ldots,m$ 
\[\pi_2(\spt_K(\Psi_i))\subset\Pi,\] so that for every choice of $(c_1,\ldots,c_m)\in\CC^m-\{\und{0}\}$, we 
also have:
\begin{equation}\pi_2(\spt_K(Q))\subset\Pi,\label{eq:linear_combination}\end{equation} where $Q$ 
is the function in (\ref{eq:rel_dep_lin}).

\begin{Lemme}
There exists exactly one partition of the set of indices ${\cal J}=\{1,\ldots,m\}$ in non-empty subsets
${\cal J}_1,\ldots,{\cal J}_p$, there exist elements
$\und{\tau}_1,\ldots,\und{\tau}_p\in{\goth H}_n^+$ pairwise $K$-linearly independent, and 
there exist elements $\varsigma_1,\ldots,\varsigma_m\in K\setminus\{0\}$ with $\varsigma_i>0$ for all $i$,
satisfying the following property.

\begin{enumerate}
\item If $i\in{\cal J}_h$, then $\und{\nu}_i=\varsigma_i\und{\tau}_h$.
\item If $i\in{\cal J}_h,j\in{\cal J}_k$ and $h\not=k$, then the supports of $\Psi_i(\und{V})$
and $\Psi_j(\und{V})$ are disjoint subsets of ${\goth H}_n^+$.
\end{enumerate}
\label{lemme:tau}\end{Lemme}

\noindent {\bf Proof.} By using (\ref{eq:support_Psi_K}), we see that $\spt_K(\Psi_i)\cap\spt_K(\Psi_j)\not=\emptyset$
implies that $\und{\nu}_i,\und{\nu}_j$ are $K$-linearly dependent.
We define the partition of the lemma from the equivalence relation induced by pairwise
$K$-linear dependence; hence (2) holds.

Thanks to (\ref{eq:support_Psi_i}), the existence of $\und{\tau}_1,\ldots,\und{\tau}_p\in{\goth H}_n^+$
such that (1) holds is guaranteed. It remaines to prove that $\varsigma_i>0$, but 
\[\varsigma_i=\frac{\nu_{i,1}}{\tau_{h,1}},\] $\nu_{i,1},\tau_{h,1}>0$ and so is $\varsigma_i$.

\medskip

We can choose two positive rational integers $q',q''$ such that, for all $i=1,\ldots,m$:
\begin{eqnarray*}
\beta_i:=q'\varsigma_i&\in &S\setminus\{0\},\\
\und{\gamma}_h:=q''\und{\tau}_h & \in &  (S\setminus\{0\})^n.
\end{eqnarray*}
Let $q=q'q''$. We have:
\begin{equation}
q\und{\nu}_i=(q'\varsigma_i)(q''\und{\tau}_h)=\beta_i\und{\gamma}_h.\label{eq:qq1q2}\end{equation}

We see that if $1\leq i\leq m_0$ then $\beta_i\in S_+$, and if $m_0+1\leq i\leq m$ then $\beta_i\in
S_\pm$. In effect, since $\tau_{h,1}'>0$ for all $h=1,\ldots,p$, the sign
of $\varsigma_i'$ (that is, the sign of $\beta_i'$) is the sign of $\nu_{i,1}'$ which is 
positive if $i=1,\ldots,m_0$ and negative if $i=m_0+1,\ldots,m$.

Let us denote:
\begin{eqnarray*}
\Upsilon_i^+(\und{u})&=&f^+(\und{\zeta}_i\und{u}^{\beta_i}),\quad\mbox{ if }i=1,\ldots,m_0,\\
\Upsilon_i^+(\und{u})&=&f(\und{\zeta}_i\und{u}^{\beta_i}),\quad\mbox{ if }i=m_0+1,\ldots,m.
\end{eqnarray*}
Using (\ref{eq:qq1q2}) we have, for $i=1,\ldots,m_0$ and $i\in{\cal J}_h$, that: 
\begin{eqnarray*}
\Psi_i(\und{v}_1^q,\ldots,\und{v}_n^q) & = & f^+(\und{\zeta}_i\und{v}_1^{q\nu_{i,1}}\cdots\und{v}_n^{q\nu_{i,n}})\\
& = & f^+(\und{\zeta}_i(\und{v}_1^{\gamma_{h,1}}\cdots\und{v}_n^{\gamma_{h,n}})^{\beta_i})\\
& = & \Upsilon_i^+(\und{v}_1^{\gamma_{h,1}}\cdots\und{v}_n^{\gamma_{h,n}}),
\end{eqnarray*}
The same computation for $i=m_0+1,\ldots,m$ implies that, for all $i=1,\ldots,m$:
\begin{equation}
\Psi_i(\und{v}_1^q,\ldots,\und{v}_n^q)=\Upsilon_i^+(\und{v}_1^{\gamma_{i,1}}\cdots\und{v}_n^{\gamma_{i,n}}).
\label{eq:identity_Upsilon}\end{equation}

We prove the lemma:

\begin{Lemme} If (\ref{eq:rel_dep_lin}) holds for non-zero complex numbers 
$c_1,\ldots,c_m$ with $Q(\und{V})$ a rational function, then
for all $h=1,\ldots,p$ there exists a rational function $Q_h(\und{u})$ such that:
\begin{equation}
\sum_{i\in{\cal J}_h}c_i\Upsilon_i^+(\und{u})=Q_h^+(\und{u}).
\label{eq:condition_refined}\end{equation}\label{lemme:from_n_to_1}
\end{Lemme} 
{\bf Proof.} 
There is no loss of generality to suppose that $c_i\not=0$ for $i=1,\ldots,m$. 

We argue by contradiction. Let us suppose that (\ref{eq:rel_dep_lin}) holds
with $Q(\und{V})$ a rational function, but there exists $1\leq h\leq p$ such that
the function $Q_h^+(\und{u})$ of (\ref{eq:condition_refined}) is irrational, i. e. not in $\CC(\und{u})$.

Since all the series $\Upsilon_i^+(\und{u})$ satisfy (\ref{eq:Upsilon_plus})
with $N_i=\beta_i^{-1}M$, lemma \ref{lemme:rank_one} applies to the linear form
(\ref{eq:condition_refined}), and from the point (2) of this lemma we see that
there exists a sequence on non-zero elements:
\[(x_i')_{i\in\NN}\subset\pi_2(\spt_K(Q_h^+))\]
such that $\lim_{i\rightarrow\infty}x_i'=0$.

The $K$-supports of 
the series \[Q_h^+(\und{v}_1^{\gamma_{h,1}}\cdots\und{v}_n^{\gamma_{h,n}})\]
are disjoint for $h=1,\ldots,p$ (lemma \ref{lemme:tau}), so that, applying (\ref{eq:identity_Upsilon}), 
there exists a sequence of points of $(\sigma_2(K)\setminus\{0\})^n$:
\[((x_{1,s}',\ldots,x_{n,s}'))_{s\in N}\subset\pi_2(\spt_K(Q(\und{V}^q)))\]
such that (\ref{eq:y1yn}) holds. Thanks to (\ref{eq:linear_combination}), lemma \ref{lemme:principle} applies, and 
$Q(\und{V}^q)$ is irrational. 
This implies that $Q(\und{V})$ is also
irrational: a contradiction.

\subsubsection{Linear forms in rational functions.\label{section:rational}}

Putting the lemmata
\ref{lemme:rank_one} and \ref{lemme:from_n_to_1} together and using lemma
\ref{proposition:equivalence_conditions}, we see that 
the validity of (\ref{eq:rel_dep_lin}) for complex numbers 
$c_1,\ldots,c_m$ not all equal to zero, with $Q(\und{V})$ a rational function, implies that for all
$1\leq h\leq p$:
\begin{equation}
\sum_{i\in{\cal
J}_h}c_i\FF_{N_i}\left(\Phi_0\left(\frac{\alpha_i}{\beta_i},\frac{\alpha_i'}{\beta_i'}
\right)\und{U}\right)=0,\label{eq:jeih}\end{equation}
where $N_i=\beta_i^{-1}M$, and some of these relation is non-trivial. 

We now consider one of the indices $h$ such that
the correspondent relation (\ref{eq:jeih}) is non-trivial: there is no loss of generality to suppose
that $h=1$. 

We only need to prove that the projective point
\begin{equation}
(c_i)_{i\in{\cal J}_1}\in\PP_{\tilde{m}-1}(\CC)
\label{eq:point_c}\end{equation}
 (with $\tilde{m}=|{\cal J}_1|$) is defined over $\QQ$: this
implies that the $\tilde{m}$-tuple $(\und{u_i})_{i\in{\cal J}_1}$
is not semi-free, as well as the original $m$-tuple ${\cal M}$, hence
completing the proof of theorem 2.

Indeed, once this proof is performed, we do as follows. We choose the point:
\[\und{v}=\Phi_0(\tau_{1,1}\xi_1+\cdots+
\tau_{1,n}\xi_n,\tau_{1,1}'\xi_1'+\cdots+\tau_{1,n}'\xi_n')\] with $\xi_i,\xi_i'\in\CC$ 
such that
\[\und{a}_i=\Phi_0(\xi_i,\xi_i'),\quad i=1,\ldots,n\]
(as in the
proof of proposition \ref{lemme:ritt}). We see that:
\begin{eqnarray*}
\und{v}^{q''} & = & \Phi_0(q''\tau_{1,1}\xi_1+\cdots+
q''\tau_{1,n}\xi_n,q''\tau_{1,1}'\xi_1'+\cdots+q''\tau_{1,n}'\xi_n')\\
& = & \und{a}_1^{\gamma_{1,1}}\cdots\und{a}_n^{\gamma_{1,n}}\in\TT(\bar{\QQ}),
\end{eqnarray*}
so that $\und{v}\in\TT(\bar{\QQ})$. Let $i\in{\cal J}_1$.
Moreover:
\begin{eqnarray*}
\und{\zeta}_i\und{v}^{q''\beta_i} & = & \und{\zeta}_i\und{a}_1^{\nu_{i,1}}\cdots\und{a}_n^{\nu_{i,n}}\\
& = & \und{\zeta}_i\und{u}_i^{\eta_i},
\end{eqnarray*}
by (\ref{eq:qq1q2}) and (\ref{eq:ritt}). In other words, (\ref{eq:jeih}) for $h=1$ 
and with rational coefficients implies that the formal series
\[(\FF_{\und{v}}(\Phi_0,\und{u}_j:\und{U}))_{j\in{\cal J}_1}\]
(definition \ref{defi:semi_free}) are not all identically zero, and
$\QQ$-linearly dependent: this is precisely what we want.

\medskip

We prove that the point $(c_i)$ of (\ref{eq:point_c}) is defined over $\QQ$.
We start with two elementary lemmata.
\begin{Lemme}[Gauss sums.]
Let $M_1\supset M_2$ be two complete $\ZZ$-modules of $K$, let $\nu\in M_2^*$. We have:
\begin{equation}
\sum_{\mu\in M_1/M_2}e(\tr(\mu\nu))=\left\{
\begin{array}{l} 0 \mbox{ if }\nu\not\in M_1^*, \\
\mbox{$[M_1:M_2]$ if $\nu\in M_1^*$,}\end{array}
\right.\label{eq:somme}
\end{equation}
The sum being indexed by a complete set of representatives of $M_1$ non equivalent modulo
$M_2$.\label{lemme:somme}\end{Lemme}

\noindent {\bf Proof.} This is well known.

\begin{Lemme}[Vandermonde matrices.]
Let $L>0$ be a rational integer, let us choose a numbering of
the sets $M^*/LM^*$ and $L^{-1} M / M $. Then,
the matrix ${\cal M}(L)$ below is non-singular:
\[{\cal M}(L)=(e(\tr(\mu\nu))_{{\tiny\begin{array}{l}
\nu\in L^{-1} M / M \\ \mu\in M^*/LM^*\end{array}}}.\]\label{lemme:matrice}
\end{Lemme}
\noindent {\bf Proof.}
The matrix is well defined and depends on the chosen numbering
of the classes of $ M^*/L M ^*$ and $L^{-1} M/ M$, but not on the representatives
chosen in any class. Let us choose a $\ZZ$-basis $(\mu_1,\mu_2)$ of $ M$, let
$(\nu_1^*,\nu_2^*)$ be the dual $\ZZ$-basis of $ M^*$.

We have, up to reorder the rows and the columns:
\begin{eqnarray*}
{\cal M}(L) & = & (e(\tr((a_1\mu_1+a_2\mu_2)(b_1\nu_1^*+b_2\nu_2^*))/L\})_{(a_1,a_2),(b_1,b_2)}\\
& = & (e((a_1b_1+a_2b_2)/L))_{(a_1,a_2),(b_1,b_2)},
\end{eqnarray*}
where the rows are indexed by the couples $(a_1,a_2)\in\ZZ^2$ with
$0\leq a_1,a_2\leq L-1$ and the columns are indexed by the couples $(b_1,b_2)\in\ZZ^2$
with $0\leq b_1,b_2\leq L-1$.

Thus the matrix ${\cal M}(L)$ is, up to permutations of rows and columns,
the Kronecker square of the Vandermonde matrix:
\[(e(ab/L))_{{\tiny 0\leq a,b\leq L-1}}.\]
Let $D$ be the determinant of this matrix.
We have $\det({\cal M}(L))=\pm D^{L}$; but $D$ is non-zero and the matrix ${\cal M}(L)$ is
non-singular.

\vspace{10pt}

We may suppose here, without loss of generality, that ${\cal J}_1=\{1,\ldots,m \}$.
The relation (\ref{eq:jeih}) for $h=1$ is equivalent to infinitely many 
linear relations indexed by the elements of $M^*$, and involving the coefficients
$c_i$ for $i=1,\ldots,m$. These relations are:
\begin{equation}
\sum_{{\tiny \begin{array}{c}i\in{\cal J}_1\mbox{ such that }\\
\nu\in N_i^*\end{array}}}c_ie(\tr(\alpha_i\nu/\beta_i))=0.
\label{eq:freeness_1}\end{equation}
Let \[{\cal N}=(A_{\nu,i})\] be the matrix of these relations;
its rows are in a one-to-one correspondence with 
the elements $\nu\in M^*$, its columns are indexed by $i=1,\ldots,m$, and the entries
are defined by:
\[A_{\nu,i}=\left\{
\begin{array}{l}
e(\tr(\alpha_i\nu/\beta_i))\mbox{ if $\nu\in N_i^*$,}\\
0 \mbox{ otherwise.}\end{array}\right.\]
We construct a matrix ${\cal N}(L)$ by cancelling almost all the rows of ${\cal N}$. 
Let $\ell_1,\ldots,\ell_{m }$ be non-zero rational integers such that
$\ell_i\alpha_i\in M$, put $g_i=[ M^*:\beta_i M^*]$ (index of complete $\ZZ$-modules). Let $N'>m$ 
be a rational integer and put \[L=N'\prod_{i=1}^{m }g_i\ell_i.\] Then we define:
\[{\cal N}(L)=(A_{\nu,i})_{{\tiny\begin{array}{l}
\nu\in M^*/LM^*\\ 1\leq i\leq m \end{array}}}.\]
The identities (\ref{eq:jeih}) for $h=1$ are equivalent to:
\[{\cal N}(L)\cdot{}^{{\rm t}}\und{c}={}^{{\rm t}}\und{0}\mbox{ in $\CC^{L^2}$}.\]
The equivalence is easy to prove because the matrix ${\cal N}$ is equal, 
up to a permutation of its rows, to a matrix made by an infinite column
which are copies of ${\cal N}(L)$.

Since the
matrix ${\cal M}(L)$ is non-singular (lemma \ref{lemme:matrice}), we have ${\cal N}(L)\cdot{}^{{\rm t}}\und{c}=
{}^{{\rm t}}\und{0}$, if and only if ${\cal M}(L)\cdot{\cal N}(L)\cdot{}^{{\rm t}}\und{c}={}^{{\rm t}}\und{0}$.
Let us compute explicitly ${\cal M}(L)\cdot{\cal N}(L)$. Let us write:
\[{\cal M}(L)\cdot{\cal N}(L)=(B_{i,\mu})_{{\tiny\begin{array}{l}
1\leq i\leq m \\ \mu\in L^{-1} M / M \end{array}}}.\] The lemma \ref{lemme:somme} can be applied
to get:
\[B_{i,\mu}=\left\{
\begin{array}{l}
[\beta_i M^*:L M^*]\mbox{ if }
\mu+\alpha_i/\beta_i\in\beta_i^{-1} M,\\
0\mbox{ if }\mu+\alpha_i/\beta_i\not\in\beta_i^{-1} M,\end{array}
\right.\]
because
\begin{eqnarray*}
B_{i,\mu} & = & \sum_{\beta\in M^*/LM^*}e(\tr(\beta\mu))A_{\beta,i}\\
& = & \sum_{\beta\in\beta_i M^*/L M^*}e(\tr(\beta(\mu+\alpha_i/\beta_i))).
\end{eqnarray*}
Let us observe that:
\[[\beta_i M^*:L M^*]=L^2/[ M^*:\beta_i M^*]=L^2/[\beta_i^{-1}M:M].\]

Let ${\cal P}(L)$ be the matrix:
\[{\cal P}(L)=(\chi_{i}(\mu))_{{\tiny\begin{array}{l}
1\leq i\leq m \\ \mu\in L^{-1} M/ M\end{array}}},\]
where $\chi_i:K/M\rightarrow\{0,1\}$ is the caracteristic function (of subsets of $K$, finite modulo $M$):
\[\chi_{i}(\mu)=\left\{
\begin{array}{l}
1\mbox{ if }\mu\in\beta_i^{-1}( M-\alpha_i)/M,\\
0\mbox{ if }\mu\not\in\beta_i^{-1}( M-\alpha_i)/M.\end{array}
\right.\]
Let us write $\und{b}=([\beta_1^{-1}M:M]^{-1}c_1,\ldots,[\beta_{m }^{-1}M:M]^{-1}c_n)$. 
We have proven that:
\[{\cal M}(L)\cdot{\cal N}(L)\cdot{}^{{\rm t}}\und{c}={\cal P}(L)\cdot{}^{{\rm t}}\und{b}=
{}^{{\rm t}}\und{0},\] or in an equivalent formulation, that the characteristic
functions $\chi_i$ of the $m $ sets of classes of $\beta_i^{-1}( M-\alpha_i)\subset L^{-1} M$
modulo $M$ are $\CC$-linearly dependent. But for any $m $-tuple of characteristic
functions of subsets of any finite set, $\CC$-linear dependence implies $\QQ$-linear dependence.
Thus $(c_i)_{i\in{\cal J}_1}\in\PP_{m -1}(\QQ)$, the relations
in (\ref{eq:jeih}) are defined over $\QQ$ and we have encountered
the required contradiction with the hypothesis of semi-freeness assumed
at the beginning; the proof of our theorem 2 is now complete.

\section{Appendix.}

In this appendix we prove a generalisation of theorem 2 for any irrational quadratic $w$'s
(section \ref{section:general_w}),
we give more informations about the property of semi-freeness (section
\ref{section:example}), we prove the corollaries 1, 2 (section \ref{section:consequences})
and finally, we give a more precise flavour of the linear relations which may occur between
complex numbers $f(\und{u}_i)$ with $\und{u}_i$ algebraic (section \ref{section:portrait}). Then,
we introduce the reader to other problems, more or less
related to Hecke-Mahler series.

\subsection{How to deal with general $w$'s in theorem 2.\label{section:general_w}}

We may suppose that $0<w<1$: it has an ordinary continued fraction development
\[w=\frac{1}{d_0+}\frac{1}{d_1+\cdots}\frac{1}{b_0+}\frac{1}{b_1+\cdots}=[0,d_0,d_1,\ldots,d_g,
\underbrace{b_0,b_1,\ldots,b_{2r-1}}_{\mbox{period}},b_0,b_1,\ldots]\]
for $g\geq 0,r>0$ and $d_0,\ldots,d_{g-1},b_0,\ldots,b_{2r-1}\in\ZZ_{>0}$. 
Following \cite{Masser:Hecke} pp. 210-211, we put:
\[D:=\sqm{0}{1}{1}{d_0}\cdots\sqm{0}{1}{1}{d_{g-1}}=\sqm{a}{b}{c}{d}\mbox{ and }
T:=\sqm{d_g}{1}{1}{0}\cdots\sqm{d_0}{1}{1}{0}=\sqm{x}{y}{z}{t}\] when $g>0$; otherwise,
we put $D=T=$ identity matrix.
Let $\theta\in K$ be defined by $w=\displaystyle{\frac{a\theta+b}{c\theta+d}}$:
we note that $\theta^{-1}$ has a purely periodic
ordinary continued fraction development:
\begin{equation}
\theta^{-1}=b_0+\frac{1}{b_1+}\frac{1}{b_2+\cdots}=[b_0,b_1,\ldots,b_{2r-1},b_0,b_1,\ldots].
\label{eq:continued_fraction}
\end{equation}
The theorems 3, 4 p. 80 of \cite{Perron:Kettenbruechen} say that the latter condition
is equivalent to (\ref{eq:perron}). We have:
\[f_w(\und{u})=f_\theta(T.\und{u})+R(\und{u})\] for some rational function $R\in\QQ(\und{u})$ defined over
couples of complex numbers $(u,v)$ such that $|u|<1$ and $|u||v|^w<1$.
(cf. \cite{Masser:Hecke}, p. 211 equation (3.6)). Thus, if ${\cal M}=((u_1,v_1),\ldots,(u_m,v_m))$
is an $m$-tuple of couples of algebraic numbers as in the hypotheses of the theorem 2, then
the complex numbers \[f_w(u_1,v_1),\ldots,f_w(u_m,v_m)\] are algebraically independent
over $\QQ$ if and only if the complex numbers \[f_\theta(T.(u_1,v_1)),\ldots,f_\theta(T.(u_m,v_m))\]
are algebraically independent over $\QQ$. We then have the following corollary of theorem 2.

\vspace{10pt}

\noindent {\bf Theorem 3.} {\em Let
${\cal M}=((u_1,v_1),\ldots,(u_m,v_m))$ be an $m$-tuple
of algebraic elements of $\TT$ such that
$|u_i|<1$ and $0<|u_i||v_i|^w<1$ for $i=1,\ldots,m$. Then $T.{\cal M}=
(T.(u_1,v_1),\ldots,T.(u_m,v_m))$ is 
semi-free with respect to $\Phi_0$ if and only if the complex numbers 
$f_w(u_1,v_1),\ldots,f_w(u_m,v_m)$ 
are algebraically independent over $\QQ$.}

\subsection{Other facts about the semi-freeness
condition\label{section:example}.}

In this subsection, we give some precisions about the semi-freeness condition.

The two lemmata below say that the condition of semi-freeness does not depend on the
choice of $\alpha,\und{v}$ in (\ref{eq:torsion_and_not}).

\begin{Lemme} The series $\FF_{\und{u}}(\Phi,\und{v}:\und{U})$ does not depend on the choice of 
a representative of the class of $\alpha$ modulo $M$ in (\ref{eq:torsion_and_not}).
\label{lemma:not_depend}\end{Lemme}
\noindent {\bf Proof.} Since $M(\und{U})^\nu M(\und{U})^\mu=M(\und{U})^{\mu+\nu}$, we may formally identify
\begin{eqnarray*}
\FF_N(\und{U}) & = & \sum_{\nu\in N^*}e( \nu Z+\nu' Z')\\
& = & \sum_{\nu\in N^*}e( \tr(\nu \und{Z}))\end{eqnarray*} for a couple of
unknowns $\und{Z}=(Z,Z')$ formally satisfying $\Phi(Z,Z')=(U,V)$.

If $\Phi(\alpha,\alpha')=\Phi(\tilde{\alpha},\tilde{\alpha}')$ for $\alpha,\tilde{\alpha}\in K$,
then $\alpha=\tilde{\alpha}+\lambda$ for some $\lambda\in M$ and given a non-zero element $\beta\in S(M)$,
we have $\tr(\tilde{\alpha}\nu/\beta)=\tr(\alpha\nu/\beta)$ for all $\nu\in\beta M^*$, thus
\[\FF_{\beta^{-1}M}\left(\Phi\left(\frac{\tilde{\alpha}}{\beta},\frac{\tilde{\alpha}'}{\beta'}\right)\und{U}\right)=
\FF_{\beta^{-1}M}\left(\Phi\left(\frac{\alpha}{\beta},
\frac{\alpha'}{\beta'}\right)\und{U}\right).\]

\begin{Lemme}
Let $(\und{u}_1,\ldots,\und{u}_m)$ be a $m$-tuple of $\TT^m$ whose coefficients
are of infinite order, let us suppose that 
there exists elements $\alpha_{i,j}\in K$, two elements of infinite order 
$\und{v}_1,\und{v}_2\in\TT$ and elements $\beta_{i,j}\in S-\{0\}$ ($i=1,\ldots,m,j=1,2$) such that:
\begin{equation}
\und{u}_i=\Phi(\alpha_{i,1},\alpha_{i,1}')\und{v}_1^{\beta_{i,1}}=\Phi(\alpha_{i,2},\alpha_{i,2}')\und{v}_2^{\beta_{i,2}},\quad
i=1,\ldots,m.
\label{eq:due_cosi}\end{equation}
Then we have:
\begin{equation}
\sum_{i=1}^mc_i\FF_{\und{v}_1}(\Phi,\und{u}_i:\und{U})=0\label{eq:inizio}\end{equation}
if and only if 
\begin{equation}\sum_{i=1}^mc_i\FF_{\und{v}_2}(\Phi,\und{u}_i:\und{U})=0.\label{eq:fine}\end{equation}
\label{lemme:two_groups}
\end{Lemme}
\noindent {\bf Proof.} Since $\und{v}_1,\und{v}_2$ are of infinite order, we can choose elements
$(z_1,z_1')$, $(z_2,z_2')\in\CC^2-\Sigma(K)$ such that:
\[\und{v}_1=\Phi(z_1,z_1'),\quad\und{v}_2=\Phi(z_2,z_2'),\]
Moreover, there exists $\tilde{\beta_1},\tilde{\beta}_2\in S-\{0\}$ such that
\[\und{v}_1^{\tilde{\beta}_1}\und{v}_2^{\tilde{\beta}_2}=\und{1},\]
so that there exists $\delta\in K-\{0\}$ and $\gamma\in K$ with
\[z_1=\delta z_2+\gamma,\quad z_1'=\delta'z_2'+\gamma'.\]
The conditions (\ref{eq:due_cosi}) are equivalent to the expressions in
$\TT$:
\[\und{u}_i=\Phi(\alpha_{i,1}+\beta_{i,1}z_1,
\alpha_{i,1}'+\beta_{i,1}'z_1')=\Phi(\alpha_{i,2}+\beta_{i,2}z_2,
\alpha_{i,2}'+\beta_{i,2}'z_2').\]
These settings imply:
\begin{equation}
\alpha_{i,1}+\beta_{i,1}\gamma-\alpha_{i,2}=\tau_i\in M
\mbox{ and }\beta_{i,1}\delta=\beta_{i,2}.\label{eq:settings}\end{equation}
Let us write $N_{i,j}=\beta_{i,j}^{-1}M$; the relation (\ref{eq:inizio}) implies the formal relation:
\[\sum_ic_i\FF_{N_{i,1}}
\left(\Phi\left(\frac{\alpha_{i,1}}{\beta_{i,1}},\frac{\alpha'_{i,1}}{\beta'_{i,1}}\right)\und{U}\right)=0.\]
The relation above is equivalent to the following formal linear relation (for formal variables
$\und{Z}_1=(Z_1,Z_1'),\und{Z}_2=(Z_2,Z_2')$ such that $Z_1=\delta Z_2+\gamma,Z_1'=\delta'Z_2'+\gamma'$):
\begin{eqnarray*}
0 & = & \sum_{i=1}^mc_i\sum_{\nu\in
N_{i,1}^*}e\left(\tr\left(\nu\left(\frac{\alpha_{i,1}}{\beta_{i,1}}+\und{Z}_1\right)\right)\right)\\ & = &
\sum_{i=1}^mc_i\sum_{\nu\in N_{i,1}^*}e\left(\tr\left(\nu\left(\frac{\alpha_{i,1}}{\beta_{i,1}}+\delta
\und{Z}_2+\gamma\right)\right)\right)\\ & = & \sum_{i=1}^mc_i\sum_{\nu\in
N_{i,2}^*}e\left(\tr\left(\nu\left(\frac{\alpha_{i,1}+\beta_{i,1}\gamma}{\beta_{i,1}\delta}+ \und{Z}_2\right)\right)\right)\\
& = & \sum_{i=1}^mc_i\sum_{\nu\in N_{i,2}^*}e\left(\tr\left(\nu\left(\frac{\alpha_{i,2}+\tau_i}{\beta_{i,2}}+
\und{Z}_2\right)\right)\right)\\
& = & \sum_{i=1}^mc_i\sum_{\nu\in N_{i,2}^*}e\left(\tr\left(\nu\left(\frac{\alpha_{i,2}}{\beta_{i,2}}+
\und{Z}_2\right)\right)\right),
\end{eqnarray*}
thanks to (\ref{eq:settings}) and because $\tau_i/\beta_{i,2}\in N_{i,2}$. Thus (\ref{eq:fine}) holds.

\begin{Lemme}
The action (\ref{eq:action_S}) only depends on the ratio $B_0/B_{1}$.
\label{lemme:ratio}\end{Lemme}
\noindent {\bf Proof.}
Let $M^\sharp=\nu M$ for some $\nu\in K-\{0\}$, let us denote $B_0^\sharp=\nu B_0,
B_1^\sharp=\nu B_1$, let $\Phi^\sharp$ be the exponential function associated to
the complete $\ZZ$-module $M^\sharp$ with the basis $(B_0^\sharp,B_1^\sharp)$.
Then the action of $S$ induced by $\Phi$ is equal to the action of $S$ induced
by $\Phi^\sharp$. 

\medskip

The lemma \ref{lemme:ratio} says that, to study the semi-freeness
condition with respect to an exponential function $\Phi$, there is no
restriction to consider complete $\ZZ$-modules $M$ with a basis $(B_0,B_1)$
such that $B_1=1$. 

\subsubsection{An explicit example.}

As we said earlier, the action of $S$ varies with the choice of the exponential function
$\Phi$ associated to it. It may happen that an $m$-tuple ${\cal M}$ is semi-free with respect
to an exponential function $\Phi$, but not semi-free with respect to another exponential function $\Psi$,
even if the underlying complete $\ZZ$-module is the same.
We explain with an example this phenomenon.

\medskip

Let us choose $K=\QQ(\sqrt{5})$ and $M=\ZZ+\epsilon\ZZ$ with $\epsilon=(1+\sqrt{5})/2$,
so that $S=M$.
We choose two bases of $M$: $(B_0,B_1)=(\epsilon,1)$ and $(C_0,C_1)=(1,\epsilon)$, so that,
with ${\goth B},{\goth C}$ associated to these bases $(B_0,B_1),(C_0,C_1)$ as in (\ref{eq:goth_B}), 
setting $T=\sqm{0}{1}{1}{0}$, we have:
\[{\goth C}=T\cdot{\goth B}.\] An easy computation shows that:
\[L={\goth C}^{-1}{\goth B}=\sqm{a}{b}{b'}{a'},\]
with $a=2B_0^*,b=-B_1^*{}'$. Let $\Phi,\Psi$ be the exponential functions
associated to ${\goth B},{\goth C}$ as in (\ref{eq:exponentielle}). 

Let $\beta\in S$ be an irrational element (for example, a unit of $S$ not equal to $\pm 1$), 
let $\und{v}=\Phi(z,z')$ be
of infinite order, let \[\und{u}=\und{v}^\beta=\Phi(\beta z,\beta'z'),\] the action of $S$ being the one 
induced by $\Phi$. The couple:
\[{\cal M}=(\und{u},\und{v})\] is not semi-free with respect to $\Phi$. Indeed, we easily check that
\[\FF_{\und{u}}(\Phi,\und{v}:\und{U})  = \FF_{\und{v}}(\Phi,\und{v}:\und{U}) = \FF(\und{U}),\] 
so that 
$\FF_{\und{u}}(\Phi,\und{v}:\und{U})$
and $\FF_{\und{v}}(\Phi,\und{v}:\und{U})$ are clearly $\QQ$-linearly dependent, regardless to the choice of
$z,z'$.

Now let $\und{w}_1=\Psi(az,a'z')$, $\und{w}_2=\Psi(bz',b'z)$, so that
\begin{eqnarray*}
\und{u} & = & \Psi({}^{\rm t}L\cdot{}^{\rm t}(z,z'))\\
& = & \und{w}_1^\beta\und{w}_2^{\beta'},\\
\und{v} & = & \und{w}_1\und{w}_2.
\end{eqnarray*}
This time, the action of $S$ is the one induced by $\Psi$ (note that $\beta'\in S$).
Since $a,b'$ are $K$-linearly independent, we may choose $z,z'$ such that the points $(az,a'z'),(bz',b'z)\in\CC^2$
are $K$-linearly independent. With this choice, we observe that 
if $\gamma,\delta\in S$ are such that
\[\und{w}_1^{\gamma}\und{w}_2^{\delta}=\und{1},\] 
then $\gamma=\delta=0$, that is, $\und{w}_1,\und{w}_2$ are multiplicatively independent with respect to
$\Psi$. In particular, if $\und{w}=\Psi(\zeta,\zeta')$ is 
a point of infinite order, the subgroup
\[\TT_{{\tiny \mbox{tors}}}\cdot\und{w}^S=\{\Phi_{{\goth
C}}(\alpha+\gamma\zeta,\alpha'+\gamma'\zeta'),\mbox{with }\alpha\in K\mbox{ and }\gamma\in S\},\]
cannot contain $\und{w}_1$ and $\und{w}_2$ at once. Thus, if $\FF_{\und{u}}(\Psi,\und{w}:\und{U})$
is non zero, then $\FF_{\und{v}}(\Psi,\und{w}:\und{U})=0$, and if
$\FF_{\und{v}}(\Psi,\und{w}:\und{U})$
is non zero, then $\FF_{\und{u}}(\Psi,\und{w}:\und{U})=0$; finally, ${\cal M}$ is semi-free
with respect to $\Psi$.

\subsection{Proof of the corollaries 1, 2.\label{section:consequences}}

We prove the corollary 1 in the case $w^{-1}=\theta$ with 
$\theta$ satisfying (\ref{eq:continued_fraction}).
If $H=\GG_m(\CC)\times\{1\}\subset\GG$ we get the theorem 1, because in this case 
$f(u,1)=f(\theta^{-1},u)$, the one-variable function of \cite{Masser:Hecke}.

The hypothesis on $H$ is equivalent to the existence of a couple of coprime
rational integers $(h,l)\in\ZZ^2-\{(0,0)\}$ such that for all $(u,v)\in H$, $u^lv^h=1$.

Let us denote $\nu=\Delta^{-1/2}(-h\theta'{}^{-1}+l)\in M^*$; we have that
for all $r\in\QQ-\ZZ$, $r\nu\not\in M^*$. Moreover, $\und{u}=\Phi_0(z,z')\in H$ (for complex numbers $z,z'$)
if and only if \begin{equation}\nu z+\nu' z'\in M.\label{eq:belongs_H}\end{equation}

Thanks to our theorem 2 we only need to prove that if ${\cal M}$ is not semi-free,
then there exists $1\leq i\not=j\leq m$ such that $\und{u}_i=\und{u}_j$. 

There is no loss of generality to suppose that there exists $\und{a}\in H(\bar{\QQ})$
of infinite order, such that:
\begin{equation}
\und{u}_i=\und{\vartheta}_i\und{a}^{\beta_i},\quad i=1,\ldots,m,\label{eq:hypo_2}
\end{equation} for torsion elements 
$\und{\vartheta}_i$ and $\beta_i\in S-\{0\}$; that is, the $S$-group generated by the $\und{u}_i$'s
has rank $1$.

\begin{Lemme} With the hypotheses above, if (\ref{eq:hypo_2}) holds,
then there exists a point of infinite order $\und{w}\in H(\bar{\QQ})$, torsion points
points $\und{\zeta_i}\in H$, and non-vanishing positive rational integers
$r_i$, such that for all $i=1,\ldots,m$:
\begin{equation}\und{u}_i=\und{\zeta}_i\und{w}^{r_i},\quad i=1,\ldots,m.\label{eq:loploppa}\end{equation}
\label{lemme:claim_1}\end{Lemme}

\noindent {\bf Proof.} We first claim that 
if two
elements $\und{u},\und{v}$ of $H(\bar{\QQ})\cap{\cal D}$ satisfy
\begin{equation}\und{u}^\beta\und{v}^\gamma=\und{1},\label{eq:beta_gamma}\end{equation}
for some $\beta,\gamma\in S-\{0\}$, 
then there exists an algebraic element $\und{w}\in
H$ of infinite order such that:
\begin{equation}
\und{u}=\und{\zeta}\und{w}^a,\quad\und{v}=\und{w}^b,\label{eq:claim_ab}\end{equation} with 
a torsion point $\und{\zeta}\in H$.

Let us consider $\xi_1,\xi_1',\xi_2,\xi_2'\in\CC$ such that:
\[\Phi_0(\xi_1,\xi_1')=\und{u},\quad
\Phi_0(\xi_2,\xi_2')=\und{v}.\]
Let $x_1,x_2,y_1,y_2,x_1',\ldots\in\RR$ be such that 
\[\xi_i=x_i+{\rm i}y_i,\quad \xi_i'=x_i'+{\rm i}y_i'.\]

The relation (\ref{eq:beta_gamma}) and (\ref{eq:belongs_H}) for $z=\xi_i,z'=\xi_i'$ imply:
\begin{eqnarray}
\beta y_1+\gamma y_2&=&\beta' y_1'+\gamma' y_2'=0,\label{eq:lo1}\\
\nu y_1+\nu' y_1'&=&\nu y_2+\nu' y_2'=0\label{eq:lo2}.
\end{eqnarray}
We see that $y_1y_1'y_2y_2'\not=0$: indeed, if for example, $y_1=0$, then by (\ref{eq:lo2})
we obtain $y_1'=0$ and this implies that $\und{u}$ lies in the boundary of ${\cal D}$, case that
has been excluded. The same argument is valid for $y_2,y_2'$.

Now, from (\ref{eq:lo1}) we see that $\beta/\gamma=-y_2/y_1\in K^\times$ and $\beta'/\gamma'=-y_2'/y_1'$.
Moreover, from (\ref{eq:lo2}), $y_2/y_1=y_2'/y_1'$, so that $\beta/\gamma=\beta'/\gamma'$; that is,
$\beta/\gamma\in\QQ^\times$.

There exist non-vanishing rational integers $p,q$ such that:
\begin{equation}
\frac{\beta}{\gamma}=\frac{q}{p}.\label{eq:pdivisoq}
\end{equation}

From the relation $\und{u}^{p\beta}\und{v}^{p\gamma}=\und{1}$
we get $\und{u}^{q\gamma}\und{v}^{p\gamma}=\und{1}$ which implies 
$\und{u}^{q}\und{v}^{p}=\und{\varphi}$ for some torsion point $\und{\varphi}\in H$.

Let $\und{w}\in H$ be any element such that $\und{w}^q=\und{v}$: it is an algebraic
point of infinite order, because so is $\und{v}$.
From $\und{u}^q=\und{\varphi}\und{w}^{-pq}$ we get $\und{u}=\und{\zeta}\und{w}^{-p}$
for some torsion point $\und{\zeta}\in H$ such that $\und{\zeta}^p=\und{\varphi}$; the claim is proved.

We end the proof of lemma \ref{lemme:claim_1} by induction on $m>0$. If
$m=1$ the property to be proved is trivially satisfied. Let us suppose that
\[\und{u}_i=\und{\varrho}_i\und{x}^{s_i},\quad i=1,\ldots,m-1.\]
We apply the claim (equality (\ref{eq:claim_ab})), to $\und{u}=\und{u}_m$ and $\und{v}=\und{x}$.
There exist a torsion point $\und{\zeta}$ of $H$ and a point of infinite order
$\und{w}\in H(\bar{\QQ})$ such that $\und{u}_m =
\und{\zeta}\und{w}^a$ and, for all
$i=1,\ldots,m-1$:
\begin{eqnarray*}
\und{u}_i &= &\und{\zeta}\und{\varrho}_i\und{w}^{as_i}\\
& = & \und{\zeta}_i\und{w}^{r_i}.
\end{eqnarray*}
For all $\und{u}\in{\cal D}$, $\und{u}\not\in{\cal D}$. Thus, the rational integers $r_1,\ldots,m$
have the same sign. Up to replace $\und{w}$ by $\und{w}^{-1}$, we obtain $r_i>0$ for $i=1,\ldots,m$.
The proof of lemma
\ref{lemme:claim_1} is complete.

\medskip

\noindent {\bf End of proof of corollary 1.} 
We apply lemma \ref{lemme:claim_1}
to the $m$-tuple \[(\und{u}_1,\ldots,\und{u}_m)\in H(\bar{\QQ})\cap{\cal D};\] let
$\und{w}$ as in (\ref{eq:loploppa}).

By hypothesis, ${\cal M}$ is not semi-free,
so that by (\ref{eq:fine}) of lemma \ref{lemme:two_groups}, there exists a non-trivial
linear dependence relation:
\[\sum_{i=1}^mc_i\FF_{\und{u}_i}(\Phi_0,\und{w}:\und{U})=0,\] with $(c_1,\ldots,c_m)\in\QQ^m-\{(0,\ldots,0)\}$.
This relation imply, as for (\ref{eq:freeness_1}), the linear relations
for all $\mu\in M^*$: 
\begin{equation}
\sum_{{\tiny \begin{array}{c}i\mbox{ such that }\\
\mu\in r_iM^*\end{array}}}c_ie(\tr(\alpha_i\mu)/r_i)=0,
\label{eq:freeness}\end{equation}
where $\alpha_1,\ldots,\alpha_m\in K$ are such that $\Phi_0(\alpha_i,\alpha_i')=\und{\zeta}_i$.

Let $\bar{N}\subset K$ be the $\ZZ$-module generated by all the elements $\alpha\in K$
such that $\tr(\alpha\nu)\in\ZZ$; then $\bar{N}=\QQ\alpha^\sharp+M$, 
for some $\alpha^\sharp\in K-\{0\}$ ($h,l$ are coprime).

We have $\alpha_1,\ldots,\alpha_m\in \bar{N}$. Let $M\subset N\subset\bar{N}$ be the complete $\ZZ$-module
generated by $M,\alpha_1,\ldots,\alpha_m$. The group $N/M$ is finite cyclic: let $\alpha\in N$
be a representative of a generator. There exists a positive rational integer $\ell$
such that $\ell\alpha\in M$ and $s\alpha\not\in M$ for $s=1,\dots,\ell-1$. There exists
$m$ rational integers $0\leq s_1,\ldots,s_m\leq\ell-1$ such that $\alpha_i\in s_i\alpha+M$.
For all $\mu\in M^*$, the relations (\ref{eq:freeness}) become:
\[\sum_{{\tiny \begin{array}{c}i\mbox{ such that }\\
\mu\in r_iM^*\end{array}}}c_ie(\tr(\alpha\mu)s_i/r_i)=0.\]
We note that the image of the map $\phi:M^*\rightarrow\QQ$ defined by $\mu\mapsto\tr(\alpha\mu)$
is $\ZZ/\ell$. Indeed, it is possible to construct a basis $(\nu,\tau)$ of $M^*$ 
such that the dual basis of $M$ is $(\kappa,\ell\alpha)$ for some $\kappa$.

Let us choose $r=r_{i_0}$ the smallest possible with $c_{i_0}\not=0$. Then, $\phi$
maps the set $rM^*-\cup_{r_j\not=r}r_jM^*$ surjectively onto $(r/\ell)\ZZ$. The relations
(\ref{eq:freeness}) become:
\[\sum_{{\tiny \begin{array}{c}i\mbox{ such that }\\ r_i=r\end{array}
}}c_ie(s_ih/\ell)=0,\mbox{ for all }h\in\ZZ.\]
The vanishing of a Vandermonde determinant implies, for two distinct indices $i,j$, that 
$r_i=r_j$ and $s_i=s_j$. 
This finally gives $\und{u}_i=\und{u}_j$, hence proving the corollary 1.

\vspace{10pt}

\noindent {\bf Proof of corollary 2.} Let us consider an $m$-tuple ${\cal M}=(\und{u}_1,\ldots,\und{u}_m)$ of algebraic
elements of $\TT\cap{\cal D}$ such that, for all $i=1,\ldots,m$, $\und{u}_i$
belongs to a given $S$-group $\Gamma$ isomorphic to $S\oplus\cdots\oplus S$, of finite 
rank $n>0$. 

Then, it is easy to see that ${\cal M}$ is semi-free if and only if for all $\und{v}
\in\Gamma(\bar{\QQ})$, if there exists a non-empty subset ${\cal J}\subset\{1,\ldots,m\}$ such that
for $j\in{\cal J}$,
$\und{u}_j=\und{v}^{\beta_j}$ with $\beta_j\in S-\{0\}$, 
then the complete modules $\beta_jM^*$
are distinct. Thus the corollary 2 is an immediate consequence of theorem 2.

\vspace{10pt}

In the introduction we quoted a very general result of Loxton and van der Poorten
which also applies to prove some properties of algebraic independence of values of Hecke-Mahler
series $f(u,v)$ with algebraic $(u,v)\in\TT\cap{\cal D}$. The application of this result 
appears on pp. 407-408 of \cite{Loxton:Arithmetic}. The corollary 2 implies these properties and suggests further
improvements.

Let $(u_1,v_1),\ldots,(u_m,v_m)$ be couples of algebraic numbers. We will say that the $m$-tuple
$((u_1,v_1),\ldots,(u_m,v_m))$ is ${\cal B}_0$-{\em free} (or {\em free}) if the following is true.
Let $(p_1,q_1,\ldots,p_m,q_m)$ be any $2m$-tuple of rational integers and let us define
the sequences $(p_j^{(k)},q_j^{(k)})=(p_j,q_j)\cdot{\cal B}_0^k$; $j=1,\ldots,m$ and $k=0,1,\ldots$. If
\[\prod_{i=1}^{m}u_i^{p_i^{(k)}}\prod_{i=1}^mv_i^{q_i^{(k)}}=1\]
for infinitely many $k\geq 0$ then $p_1=q_1=\cdots=p_m=q_m=0$.

\vspace{10pt}

\noindent {\bf Theorem 4. (Loxton and van der Poorten.)}
{\em Let $((u_1,v_1),\ldots,(u_m,v_m))$ be an $m$-tuple of couples of algebraic 
numbers such that
$0<|u_j||v_j|^\theta<1$ and $|u_j|<1$ for all $j=1,\ldots,m$. Then if the 
$m$-tuple $((|u_1|,|v_1|),\ldots,(|u_m|,|v_m|))$ is free, the complex numbers $f(u_1,v_1),\ldots,
f(u_m,v_m)$ are algebraically independent over $\QQ$.}

\vspace{10pt}

This theorem contains a theorem of Mahler in \cite{Mahler:Arithmetische} (case $m=1$).
It is worth to remark that the proof of the lemma \ref{lemme:proprieteA}
says also that an $m$-tuple $(\und{u}_1,\ldots,\und{u}_m)$ is 
free if and only if $\und{u}_1,\ldots,\und{u}_m$ are multiplicatively 
independent. An $m$-tuple $(\und{u}_1,\ldots,\und{u}_m)$ of multiplicatively 
independent elements of $\TT(\bar{\QQ})\cap{\cal D}$ clearly 
satisfies the hypotheses of the corollary 2.

If an $m$-tuple ${\cal M}=((u_1,v_1),\ldots,(u_m,v_m))$ of couples of algebraic numbers is such that
the $m$-tuple $((|u_1|,|v_1|),\ldots,(|u_m|,|v_m|))$ is free, then ${\cal M}$ is
also free: this shows that our corollary 2 implies the theorem 4.

Note that there exist $m$-tuples $(\und{u}_1,\ldots,\und{u}_m)$ which are
semi-free but not free. An example is given by any sub-$\no(\beta)$-tuple of ${\cal M}$ as in
(\ref{eq:an_example}), so that our theorem 2 is really more general than theorem 4.

\subsection{A portrait of the relations.\label{section:portrait}}

We divide the relations in two different types.

\subsubsection{Relations of generic type.}

There are linear relations that we call Ògeneric"
which hold for all of $w\in\RR_{>0}$ at once, regardless to the quadraticity of $\theta$ and the algebraicity 
of the points. 
These inconditional relations are all
homogeneous and defined over $\QQ$. If $(u,v)$ is a couple of
complex numbers in $\TT\cap{\cal D}$, the simplest example is the relation (\ref{eq:masser}), with 
$m=5$ and \[{\cal M}=((u,v),(u,-v),(-u,v),(-u,-v),(u^2,v^2)),\] which is not semi-free
(note however, that $((u,v),(u,-v),(-u,v),(-u,-v))$ is semi-free, as well as any other 
sub-$k$-tuple of ${\cal M}$ with $1\leq k\leq 4$).

These relations all arise from the linear relations of rational functions:
\[\sum_{i=1}^m c_if_\infty(\zeta_i u^{b_i},\zeta_i v^{b_i})=0,\]
where $c_1,\ldots,c_m$ are rational numbers, $\zeta_1,\ldots,\zeta_m$ are roots of unity,
$b_1,\ldots,b_m$ are positive rational integers, and $f_\infty=\lim_{w\mapsto\infty}f_w$ 
is the rational function
\[f_\infty(u,v)=\frac{uv}{(1-u)(1-v)}.\]

\subsubsection{Relations of special type.}

There are linear relations that we call Òspecial", which hold for irrational quadratic $w$'s only.
these relations are almost all non-homogeneous, and defined over $\bar{\QQ}$ but not necessarily on $\QQ$.
If $w=\theta$ satisfies (\ref{eq:perron}), the simplest example of these relations is
provided for $m=2$ by the functional equation (\ref{eq:functional_equation}) (that is, (\ref{eq:fu_eq})). 

Since $R_\eta$ is a rational function defined over
$\QQ$ and converging on ${\cal D}$, if $\und{v}\in\TT(\bar{\QQ})\cap{\cal D}$,
the two complex numbers $f(\und{v})$ and $f(\und{v}^\eta)$ differ by the algebraic number
$R_\eta(\und{v})$. Of course \[{\cal M}=(\und{v},\und{v}^{\eta})\] is not
semi-free.

But there are many other relations: we give here a concrete example. Let us
consider an element $\beta\in S$ such that $\beta>\beta'>0$. Let $\Ker(\beta)$ be the kernel of the
isogeny $\TT\rightarrow\TT$ given by $\und{u}\mapsto\und{u}^\beta$:
this kernel has $\no(\beta)$ elements, indeed there is an isomorphism of groups $\Ker(\beta)
\cong (\beta^{-1}M/M)$.
Choose any element $\und{v}\in\TT_m(\bar{\QQ})\cap{\cal D}$ and consider the
$(\no(\beta)+1)$-tuple
of elements of $\TT(\bar{\QQ})\cap{\cal D}$:
\begin{equation}
{\cal
M}=(\und{u}_0,\und{u}_1\ldots,\und{u}_{\no(\beta)})=(\und{v}^\beta,\und{\zeta}_1\und{v},\ldots,
\und{\zeta}_{\no(\beta)}\und{v}),
\label{eq:an_example}\end{equation}
where $\und{\zeta}_1,\ldots,\und{\zeta}_{\no(\beta)}$ are all the elements of $\Ker(\beta)$. Clearly
${\cal M}$ is not semi-free in this case (but the $\no(\beta)$-tuple 
$(\und{u}_1\ldots,\und{u}_{\no(\beta)})$ is semi-free).

Let $\eta\in S_+$ be any unit fixing each element of $\Ker(\beta)$, such that $\eta>1$. Then we have the linear
relation:\[\sum_{i=1}^{\no(\beta)}R_{\eta,M}(\und{\zeta}_i\und{u})=\no(\beta)R_{\eta,\beta^{-1}M}(\und{u}),\]
which is easily checked by using, for example, the lemma \ref{lemme:somme} of this text. This relation
implies a relation of special type:
\[\no(\beta)f(\und{u}_0)-\sum_{i=1}^mf(\und{u}_i)=\lambda\in\bar{\QQ}.\] Note that $\lambda$ might be non-zero
because $\beta$ is chosen to be irrational, and it is easy to produce examples with $\lambda\not=0$
in this case. 

The relations of special type all arise from the rational linear relations of rational
functions (\ref{eq:relation0}) and (\ref{eq:relation1}) of lemma \ref{lemme:three_conditions}
 when $\beta_1,\ldots,\beta_m$ are elements of $S$ not all rational,
and $\eta>1$ a suitable unit of $S_+$.

In this text we have studied the case of $w$ quadratic only, essentially because we need the functional equation
(\ref{eq:functional_equation}) to apply the classical Mahler method. 
This leads to many other relations of special type. 

For an irrational non-quadratic $w>0$, only relations
of generic type should hold, and we could expect that if $\und{u}_1,\ldots\und{u}_m$ are couples of
algebraic numbers in $\TT$ close enough to $(0,0)\in\CC^2$, such that for two distinct irrationals $w_1,w_2>0$
the two relations:
\[\sum_{i=1}^mc_if_{w_j}(\und{u}_i)=0\] hold for rational numbers $c_1,\ldots,c_m$ not all zero and $j=1,2$,
then the two relations both arise from a single relation of generic type.

\subsection{Other problems.}

Here are
some other problems that we can solve.
We can use the ideas of this paper to prove that our theorem holds with $f$ replaced by
the following bi-variable Fredholm series:
\[g(\und{u})=\sum_{k=0}^\infty P(\und{u}^{\eta^k}),\]
where $P(u,v)=uv$, and we can even get a specific description of the
transcendence degree of fields such as:
\[\QQ(f(\und{u}_1),\ldots,f(\und{u}_n),g(\und{u}_1),\ldots,g(\und{u}_n)),\] for couples of
algebraic numbers $\und{u}_i$ such that the complex numbers above make sense.

Our theorem 2 can also be extended to a more general result
describing in a completely explicit way all the algebraic dependence relations over $\QQ$
of the complex numbers:
\[\frac{\partial^{l+h} f(\und{u}_i)}{\partial u^l\partial v^h},\]
for $h,l\geq 0$. This result would also contain a theorem of Nishioka (cf. theorem 3.4.8
on p. 102 of \cite{Nishioka:Mahler}).

We can say something about the problem of the
algebraic independence of complex values at algebraic points of general transcendental
solutions $f$ of functional equations such as (\ref{eq:functional_equation}), with $R_\eta$
replaced by a rational function in $n$ complex variables. 
 
We can consequently refine certain results of Loxton and van der Poorten
(for example those on pp. 407-408 of \cite{Loxton:Arithmetic}) in the following case:
the matrix ${\cal B}\in$Mat${}_n(\NN)$ has all of its eigenvalues $\lambda_1,\ldots,\lambda_n$ conjuguate
over $\QQ$ and moreover the eigenvalues $\lambda_1,\lambda_2$ are
real such that $0<\lambda_1<1<\lambda_2$, and the other
eigenvalues 
$\lambda_j$ are complex with $\lambda_1<|\lambda_j|<\lambda_2$ for $j=3,\ldots,n$.

\vspace{10pt}

As a conclusion of this text, here are some problem we cannot solve yet.

\vspace{10pt}

\noindent {\bf Problem (1).} 
Find the analogue of theorem 2 for the Hecke geometric series $G(\und{u})=\lim_{k\mapsto -\infty}
f(\und{u}^{\eta^{k}})$. This series is not locally analitic, thus the proposition
\ref{propo:loxton} cannot be applied. However, the techniques introduced in \cite{Corvaja:Application}
may successfully be applied to get at least transcendence results of complex values of these
functions at couples of algebraic numbers. 

\vspace{10pt}

\noindent {\bf Problem (2).}
Extend theorem 2 to the natural analogue of the series $f$ obtained by replacing
the quadratic number field $K$ by any totally real number field of degree $d>2$.
These are series in $d$ complex variables. These series satisfy much more complicated systems
of functional equations, as soon as $d>2$, and the proposition \ref{propo:loxton} cannot be applied to this situation.
Nevertheless, the techniques of \cite{Corvaja:Application} allow once more to reach transcendence results.

\vspace{10pt}

\noindent {\bf Problem (3).}
Find an analogue of the corollary to the theorem 2 
with $H$ replaced by a connected algebraic curve contained in $\TT$
and defined over $\bar{\QQ}$.

\begin{small}\end{small}
\noindent{\small Federico Pellarin,\\
\noindent L.M.N.O., UniversitŽ de Caen,\\
\indent Campus II - Boulevard MarŽchal Juin,\\
\indent BP 5186 - F14032 Caen Cedex.\\
\end{document}